\documentclass[11pt,twoside]{amsart}
\usepackage[dvipsnames]{xcolor}
\usepackage{hyperref}
\hypersetup{colorlinks=true, linkcolor=Blue, citecolor=Maroon}
\usepackage{amsthm, amsmath, amscd, amssymb,centernot,txfonts}
\usepackage{tikz-cd}
\usepackage{cancel}
\usepackage{lipsum}
\usepackage{multicol}
\usepackage{float}
\usepackage[normalem]{ulem}
\usepackage{amsfonts}
\usepackage{mathrsfs}
\usepackage{parskip}
\usepackage[all]{xy}
\usepackage[left=2.5cm,top=2.5cm,bottom=3cm,right=2.5cm]{geometry}
\usepackage{dirtytalk}
\usepackage{verbatim}
\usepackage{mathtools}
\usetikzlibrary{shapes,arrows}
\usepackage{etoolbox}
\usepackage{dynkin-diagrams}
\usepackage{braket}
\usepackage{subfloat}
\usepackage{enumitem}
\usepackage{chngpage}
\usepackage{adjustbox}
\usepackage[utf8]{inputenc}
\usepackage{fourier} 
\usepackage{array}
\usepackage{makecell}
\usepackage{subcaption} 
\usepackage{caption} 
\usepackage[linktocpage=true]{}

\usepackage{color, colortbl}
\definecolor{LightCyan}{rgb}{0.88,1,1}

\usepackage[normalem]{ulem}

\usepackage[first=0,last=9]{lcg}

\newcommand{\stkout}[1]{\ifmmode\text{\sout{\ensuremath{#1}}}\else\sout{#1}\fi}

\theoremstyle{plain}
\numberwithin{equation}{section}
\newtheorem{theorem}{Theorem}[section]
\newtheorem{proposition}[theorem]{Proposition}
\newtheorem{lemma}[theorem]{Lemma}
\newtheorem{corollary}[theorem]{Corollary}

\newtheorem{set-up}[theorem]{Set-up}

\theoremstyle{definition}
\newtheorem{remark}[theorem]{Remark}

\newtheorem{example}[theorem]{Example}
\newtheorem{definition}[theorem]{Definition}

\newcommand*{\QEDB}{\hfill\ensuremath{\square}}

\tikzstyle{decision} = [diamond, draw, , 
    text width=4.5em, text badly centered, node distance=3cm, inner sep=0pt]
\tikzstyle{block} = [rectangle, draw, , 
    text width=10em, text centered, rounded corners, minimum height=2em]
\tikzstyle{block1} = [rectangle, draw, , 
    text width=5em, text centered, rounded corners, minimum height=2em]    
\tikzstyle{line} = [draw, -latex']
\tikzstyle{cloud} = [draw, ellipse,, node distance=3cm,
    minimum height=2em]

\begin{document}

\title[Tautological families of cyclic covers of projective spaces and rationality of stacks]{Tautological families of cyclic covers of projective spaces }

\author[P. Kundu]{Promit Kundu}
\address{Institute of Mathematical Sciences, ShanghaiTech University, Shanghai, China}
\email{kundupromit63@gmail.com}

\author[J. Mukherjee ]{Jayan Mukherjee} \thanks{Corresponding author: Jayan Mukherjee, email id: mukherjeejayan@gmail.com}
\address{Department of Mathematics, Oklahoma State University, Stillwater, USA}
\email{jayan.mukherjee@okstate.edu}

\author[D. Raychaudhury]{Debaditya Raychaudhury}
\address{Department of Mathematics, University of Arizona, Tucson, USA}
\email{draychaudhury@math.arizona.edu}

\subjclass[2020]{14A20, 14C15, 14D22, 14D23, 14E08, 14F22}
\keywords{Tautological families, cyclic covers, rationality, linearization, Brauer--Severi schemes}

\maketitle
\begin{abstract}
 In this article we study the existence of tautological families on Zariski open subsets of the coarse moduli space of simple cyclic covers of projective spaces of arbitrary dimension and the coarse moduli space of Galois triple covers of the projective line. Both are natural generalizations of the study of such families of hyperelliptic curves carried out in \cite{GV1} and \cite{GV2}.  We completely classify the existence or non-existence of such families, give a measure of how large the open set could be and also give a sufficient criterion for their uniqueness. As a by-product of our study, we classify exactly when the Brauer--Severi scheme associated to the cyclic cover is the projectivization of a vector bundle over the base. Since the existence of a tautological family implies the neutrality of the generic gerbe, as a corollary we determine when the moduli stack of cyclic covers of $\mathbb{P}^2$ or $\mathbb{P}^1$ is rational or not, i.e., when the stack is birational to a weighted projective stack or more generally a toric Deligne-Mumford stack (see \cite{BH08}). 
\end{abstract}

\section{Introduction}

The question of existence of tautological families is a fundamental  \textcolor{black}{question} for moduli problems. For example, it is a starting point of several conjectures, such as the Franchetta conjecture on moduli of curves (now a theorem of Harer-Arbarello-Cornalba, see \cite{AC87}, \cite{Har83}), $K3$ surfaces, and hyperk\"ahler manifolds. For moduli functors parametrizing objects, such that the general object has trivial automorphism group, the existence of tautological families on open sets of the coarse moduli is a consequence of representability of an open subfunctor. This family therefore uniquely exists on the entire open subspace of the coarse moduli space parametrizing objects without automorphisms. An example of such is $M_g$, the moduli space of smooth curves of genus $g$ with $g \geq 3$ and the open set is denoted by $M_g^0$. \par
 We recall that a general elliptic curve has an automorphism of order two. In \cite{Mum63}, Mumford constructed explicitly, a tautological family over $M_{1,1}^0$, the open subset of $M_{1,1}$, consisting of elliptic curves without extra automorphisms. In the case of $H_g$, the moduli of hyperelliptic curves of genus $g$, a general hyperelliptic curve again has an automorphism of order two corresponding to the hyperelliptic involution. If $\textrm{char}(\mathbb{k}) \neq 2$, in \cite{GV1}, Gorchinskiy and Viviani show that $H_g$ has a tautological family on an open set if and only if $g$ is odd. They further show that such a family, even when it exists cannot be extended to the open subscheme $H_g^0$, of hyperelliptic curves without extra automorphisms contrary to the case of either $M_g$, $g \geq 3$, or $M_{1,1}$.
It was shown by Ran in \cite{Ran91} that there is no tautological family on the entire moduli of framed hyperelliptic curves, (i.e, hyperelliptic curves $C \to \mathbb{P}^1$ with a fixed $2:1$ map onto $\mathbb{P}^1$), but one such exists after removing a hypersurface of odd degree.  \par
It has been shown by Markman in \cite{Mar21} that tautological families exist over every irreducible component of the moduli space of marked irreducible holomorphic symplectic manifolds. \par
In \cite{BH08}, Biswas and Hoffmann prove the existence of tautological families (or Poincar\'e families) over the moduli of certain symplectic bundles on a curve. \par
Existence of tautological families is an integral part of the concept of rationality of a Deligne-Mumford stack. A DM stack is rational if and only if its coarse moduli is rational and there exists a tautological family over an open subset of its coarse moduli space. If the generic automorphism group of a DM stack is abelian, it is rational if and only if it is birational to a Toric Deligne-Mumford stack (see ~\ref{rationality of a stack}, (ii)). In \cite{GV1}, using the fact that the moduli space $H_g$ is rational (see \cite{Bog}), the authors conclude that the stack of hyperelliptic curves of genus $g$ is rational if and only if $g$ is odd. In \cite{BH08}, using the existence of the tautological families, it was shown that certain moduli spaces of symplectic bundles on curves are rational. The definition of rationality of stacks has its roots in the articles \cite{KS99} and \cite{H10} where a more general notion of birationally linear map of stacks was introduced in the context of the fundamental question of rationality of coarse moduli of vector bundles on a curve with fixed determinant when rank and degree are co-prime.

\vspace{0.2cm}

\noindent\textbf{Objective of the article.} The objective of this article is two-fold. We determine the existence/non-existence of tautological families on open subsets of 
\begin{itemize}
    \item[(i)] $M_{n,r,d}:=$ the coarse moduli of smooth simple cyclic covers of degree $r$ and branch degree $rd$ over a projective space of dimension $n$, i.e., the coarse moduli space of the stack $\mathscr{H}_{n,r,d}$.
    \item[(ii)] $M_{1,3,d_1,d_2}:=$ the coarse moduli of smooth cyclic triple covers on $\mathbb{P}^1$ with branch divisor $D_1+D_2$, where $D_1$ is of degree $2d_1-d_2$ and $D_2$ is of degree $2d_2-d_1$, which is the coarse moduli of a quotient stack $\mathscr{H}_{1,3,d_1,d_2}$. 
\end{itemize}
Both of these stacks are compelling generalizations of the stack of hyperelliptic curves, the first one in terms of dimension of the base and the degree of the cover and the second one in terms of the complexity of the Galois cover. Note that a cyclic triple cover is not necessarily simple cyclic.

\vspace{0.2cm}

\noindent\textbf{Statements of the main results.} We state our results for the stack of smooth simple cyclic covers. In the pioneering paper \cite{AV}, Arsie and Vistoli give a quotient stack structure to the category $\mathscr{H}_{n,r,d}$ fibred over $(\textrm{Sch}/\mathbb{k})$ by groupoids consisting of simple cyclic covers of degree $r$ and branch degree $rd$, over the projective space of dimension $n$ and isomorphisms between them. 

\smallskip

Before we state our main result, we recall the definition of a Brauer--Severi scheme.
\begin{definition}
A {\bf Brauer–Severi scheme} of dimension $n$ over a scheme $S$ is a proper smooth morphism $P \to S$ such that there exists an \'etale cover $S'\to S$ such that the base change $P\times_{S}S'$ is isomorphic to $\mathbb{P}_{S'}^n$ (see \cite{OM}).

\end{definition}

Note that by \cite{GroI}, \cite{GroII}, A scheme $P \to S$ is a Brauer–Severi scheme if and only if $P \to S$
is a finitely-presented proper flat morphism and if all of
its geometric fibres are isomorphic to projective spaces. 

\smallskip

Our main result is the following in which we also compute the Picard group of the coarse moduli $M_{n,r,d}$.

\begin{theorem}\label{main1intro}
($=$ Theorem \ref{tautological family}, Proposition \ref{tautological family in any dimension}, Theorem \ref{Picard Group of the coarse moduli},) Let $M_{n,r,d}^0$ be the open subscheme of $M_{n,r,d}$ parametrizing simple cyclic covers of $\mathbb{P}^n$, whose automorphism group is $\mu_r$. Then the following statements hold.
\begin{itemize}
    \item[(1)] 
    \begin{itemize}
    \item[(i)] There exists a tautological family of cyclic covers over an open set $U$ of the coarse moduli $M_{n,r,d}$ if and only if $\textrm{gcd}(rd, n+1) \mid d$.
    \item[(ii)] When such a family exists, the Brauer--Severi scheme $P \to U$ associated to the tautological family is trivial if and only if $\textrm{gcd}(rd, n+1) = 1$.
    \item[(iii)] \textcolor{black}{Assume that the codimension of the complement of $M_{n,r,d}^0$ in $M_{n,r,d}$ is greater than or equal to $2$. Then} the family when it exists is unique when restricted to the complement of all hyperplane sections in $M_{n,r,d}^0$ in the linear systems of each unique line bundle of order $k$ where $k \mid \textrm{gcd}(r,n+1)$. In particular, if $\textrm{gcd}(r,n+1) = 1$, then for those open sets over which a tautological family exists, such a family is unique.
\end{itemize}
    
    \item[(2)] \textcolor{black}{There does not exist a tautological family on $M_{n,r,d}^0$}
   
     \item[(3)] If $rd \geq 4$, $\textrm{Pic}(M_{n,r,d}) = 0$.
\end{itemize}
\end{theorem}

Since the existence of a tautological family implies the neutrality of the generic gerbe, we have, by combining Theorem ~\ref{main1intro} with known results on rationality of moduli of plane curves (see \cite{BB10}, \cite{BBK09}), the following corollary:

\begin{corollary} ($=$ Corollaries \ref{rath1}, \ref{rath2})
    \begin{itemize}
    \item[(i)]  if $rd \geq 4$, the stack $\mathscr{H}_{1,r,d}$ is rational if and only if either $d$ is even or $rd$ is odd,
    \item[(ii)]
    \begin{itemize}
    \item[(a)] if $rd \geq 49$ the stack $\mathscr{H}_{2,r,d}$ is rational if and only if either $3\mid d$ or $3\nmid rd$. 
    \item[(b)] If $(r,rd) \in S:=\left\{(3,6), (3,12), (3,15), (9,18), (3,24), (6,24), (12,24), (3,48), (6,48)\right\}$, then $\mathscr{H}_{2,r,d}$ is not rational.
\end{itemize}
\end{itemize}
\end{corollary}

We now state our results for the stack of smooth cyclic (not necessarily simple cyclic) triple covers over $\mathbb{P}^1$. Let as before $\mathscr{H}_{1,3,d_1,d_2}$ denote the category fibred over the category $(\textrm{Sch}/\mathbb{k})$ by groupoids consisting of cyclic triple covers over $\mathbb{P}^1$ branched along divisors $D_1$ and $D_2$ of degree $2d_1-d_2$ and $2d_2-d_1$ respectively and isomorphisms between them. 
\begin{theorem}\label{main2intro}
($=$ Theorem \ref{tautgal}, \textcolor{black}{Proposition} ~\ref{non-existence of global tautological family for cyclic}, Remark ~\ref{zlt13}, Proposition \ref{M_{1,3,d_1,d_2}}) Let $M_{1,3,d_1,d_2}^0$ be the open subscheme of $M_{1,3,d_1,d_2}$ parametrizing cyclic triple covers of $\mathbb{P}^1$, whose automorphism group is $\mu_3$. Assume $2d_1-d_2 \geq 4$ and $2d_2-d_1 \geq 4$. Then the following statements hold.
\begin{itemize} 
 \item[(1)] 
 \begin{itemize}
    \item[(i)] There exists a tautological family of cyclic triple covers over some open subscheme of the coarse moduli $M_{1,3,d_1,d_2}$ of $\mathscr{H}_{1,3,d_{1},d_{2}}$.
    \item[(ii)] The family is unique when restricted to the complement of all hyperplane sections in $M_{1,3,d_{1},d_{2}}^0$ in the linear systems of each line bundle of order 3.
    \item[(iii)] In particular, if one of the following holds:
    \begin{itemize}
        \item[(i)] $2 \mid  \textrm{gcd}(d_{1},d_{2})$ and $3 \mid \textrm{gcd}(l_{1}-2,l_{2}-2)$; or
        \item[(ii)] \textcolor{black}{$\textrm{gcd}(d_{2},2)=1$ and $$3(2(l_2(l_1+1)+(l_2-2)(d_2+1))-(2l_2-1)(2d_1+1))\nmid$$ $$\text{gcd}\left(\begin{array}{l} (l_2-1)(4d_2-5d_1(d_2+1)-4d_1^{2})+(l_2-1
)(d_1+2)(4d_1-5d_2),\\ -2(l_2-1)(l_1-1)(2d_1+1)+4(l_1-1)(l_2-1)(d_1+2),\\
2(l_1-1)(4d_2-5d_1(d_2+1)-4d_1^{2})+(l_1-1)(2d_1+1)(-5d_2+4d_1)\end{array}\right),$$}
    \end{itemize} 
    then for those open sets over which a tautological family exists, such a family is unique.
\end{itemize}
  \item[(2)] There does not exist a tautological family on the open subscheme $M_{1,3,d_1,d_2}^0$ of $M_{1,3,d_1,d_2}$
  \item[(3)] $\mathscr{H}_{1,3,d_{1},d_{2}}$ is unirational. $M_{1,3,d_1,d_2}$ is unirational and is fibred over a rational base by fibres birational to a homogeneous space. If $\textrm{char}(\mathbb{k}) = 0$, $M_{1,3,d_1,d_2}$ is rationally fibred over a rational base.
\end{itemize}
\end{theorem}


In view of our results, it would be interesting to study Franchetta type results on the tautological families thus obtained (for example see \cite{OG13}). It will also be interesting to know the divisibility of the modular map of a family of maximum variation of moduli with a section for either $\mathscr{H}_{n,r,d}$ or $\mathscr{H}_{1,3,d_1,d_2}$ when the modular map is generically finite.

\vspace{5pt}

\noindent\textbf{Notation and conventions.} Throughout this article, we work over an algebraically closed base field $\mathbb{k}$ of characteristic zero or $p>0$ which:
\begin{itemize}
\item[(1)] in Section \ref{2} satisfies $\textrm{gcd}(p,2rd)=1$ and $r\geq 2$,
\item[(2)] in Section \ref{3} satisfies $\textrm{gcd}(p,2l_{1}l_{2})=1,$ where $l_{1}=2d_{1}-d_{2}$ and $l_{2}=2d_{2}-d_{1}.$
\end{itemize}
Recall that a homogeneous polynomial $F(x_0,..,x_n)$ in variables $x_0,..,x_n$ is called {\bf smooth} if given any point $P \in \mathbb{A}^{n+1}-\{0\}$, $\exists j$ such that $\frac{\partial F}{\partial x_j}(P) \neq 0$. Further, $\mathbb{A}_{sm}(n,k)$ denotes the quasi-affine scheme consisting of smooth homogeneous polynomials of degree $k$ in $n+1$ variables. 

\vspace{5pt}

\noindent\textbf{Acknowledgements.} We thank Professor Purnaprajna Bangere, Professor Indranil Biswas, Professor Patricio Gallardo, Professor Brendan Hassett, Professor Yunfeng Jiang, Professor D.S. Nagaraj, and Professor Angelo Vistoli for valuable comments and suggestions. The second author was supported by the National Science Foundation, Grant No. DMS-1929284 while in residence at ICERM 
in Providence, RI, as part of the ICERM Bridge program. The research of the third author was supported by a Simons Postdoctoral Fellowship from the Fields Institute for Research in Mathematical Sciences. All authors declare that they have no conflicts of interest.

\section{Stacks of simple cyclic covers of \texorpdfstring{$\mathbb{P}^n$}{}}\label{2}

Throughout this section, we work with the stack of simple cyclic covers of projective spaces with the objective of proving Theorem ~\ref{main1intro}.

\subsection{Rigidification of the stack of cyclic covers over a projective space}\label{sub1} 

We formally introduce the main object of our study, namely the fibered category $\mathscr{H}_{n,r,d}$. 

\begin{definition}\label{definition of the stack}
Let $\mathscr{H}_{n,r,d}$ denote the fibred category 
\color{black}
\begin{itemize}
    \item[(a)] whose objects over any $\mathbb{k}$-scheme $S$ consists of a Brauer--Severi scheme $P \rightarrow S$ of relative dimension $n$, a line bundle $\mathscr{L}$ on $P$ which restricts to degree $-d$ on a fibre of every geometric point $s$ of $S$ and an injective homomorphism $i : \mathscr{L}^{\otimes r} \rightarrow \mathscr{O}_P$ that remains injective with smooth associated divisor when restricted to any geometric fibre
    \item[(b)] and morphisms consist of Cartesian diagrams between Brauer--Severi schemes preserving the line bundle and the injective homomorphism.  
\end{itemize}
\color{black}

\end{definition}

As has been shown in \cite{AV}, $\mathscr{H}_{n,r,d}$ is isomorphic to the stack of smooth simple cyclic covers of degree $r$, over $\mathbb{P}^n$ branched along a divisor of degree $rd$ (\cite{AV}, Remark 3.3). The quotient stack structure of $\mathscr{H}_{n,r,d}$ has been provided by Arsie and Vistoli in \cite{AV}, we recall the precise result below.

\begin{theorem}\label{AV} (\cite{AV}, Corollary 4.2)
$\mathscr{H}_{n,r,d} \cong [\mathbb{A}_{sm}(n,rd) / (GL_{n+1}/\mu_d)]$ where the action of $GL_{n+1}/\mu_d$ on $\mathbb{A}_{sm}(n,rd)$ is given by $[A]\cdot f(x)=f(A^{-1}x)$.
\end{theorem}

\begin{remark}
Note that given an object $\zeta = (P \to S, \mathscr{L}, i: \mathscr{L}^{\otimes r} \to \mathscr{O}_P) \rightarrow S \in Ob(\mathscr{H}_{n,r,d})$ we have that multiplication by an $r-$th root of unity in $\mathscr{O}_S^*$ is an automorphism of $\zeta \to S$. Hence  $$\mu_r \subseteq \textrm{Aut}(\zeta \to S).$$
\end{remark}

We now construct a stack $\mathscr{D}_{n,r,d}$, which as we shall see in the next section, is a stack intermediate to $\mathscr{H}_{n,r,d}$ and its coarse moduli space with generic trivial stabilizer. This key object was used to study the stack of hyperelliptic curves in \cite{GV1} and \cite{GV2}. We shall see in Remark \ref{rigidification} that $\mathscr{D}_{n,r,d}$ is the $\mu_r$-rigidification of $\mathscr{H}_{n,r,d}$.

\begin{definition}\label{definition of the orbifold}
Let $\mathscr{D}_{n,r,d}$ denote the fibred category 
\color{black}
\begin{itemize}
    \item[(a)] whose objects over any $\mathbb{k}$-scheme $S$ consist of a Brauer--Severi scheme $P \to S$ of relative dimension $n$, a line bundle $\mathscr{L}$ which restricts to degree $-rd$ on a fibre of every geometric point $s \in S$ and an injective homomorphism $i: \mathscr{L} \to \mathscr{O}_P$, that remains injective with smooth associated divisor, when restricted to any geometric fibre
    \item[(b)]  morphisms consist of Cartesian diagrams between Brauer--Severi schemes preserving the line bundle and the injective homomorphism.
\end{itemize}
\color{black}

\end{definition}

The techniques of \cite{AV} gives a quotient structure to $\mathscr{D}_{n,r,d}$.

\begin{theorem}\label{orbifold}
The stack $\mathscr{D}_{n,r,d} \cong [\mathbb{A}_{sm}(n,rd) / (GL_{n+1}/\mu_{rd})]$ with the action given by $[A] \cdot f(x) = f(A^{-1} x)$.
\end{theorem}

\color{black}
\noindent\textbf{Proof.} Consider the auxiliary fibred category $\widetilde{\mathscr{D}}_{n,r,d}$ whose objects over a $\mathbb{k}$-scheme $S$ consist of an object $(P \to S, \mathscr{L}, i: \mathscr{L} \to \mathscr{O}_P)$ in $\mathscr{D}_{n,r,d}$ along with an isomorphism $\phi: (P, \mathscr{L}) \cong (\mathbb{P}_S^n, \mathscr{O}_{\mathbb{P}_S^n}(-rd))$. Note that there is forgetful map from $\widetilde{\mathscr{D}}_{n,r,d} \to \mathscr{D}_{n,r,d}$. We define a base preserving functor from $\widetilde{\mathscr{D}}_{n,r,d} \to \mathbb{A}_{sm}(n,r,d)$. For any object of $\mathscr{\widetilde{D}}_{n,r,d}(S)$ consider the composite homomorphism $\phi \circ i \circ \phi^{-1}: \mathscr{O}_{\mathbb{P}_S^n}(-rd)) \to \mathscr{O}_{\mathbb{P}_S^n}$ which gives an element of $\mathbb{A}_{sm}(n,rd)(S)$. Conversely given an element of $\mathbb{A}_{sm}(n,rd)(S)$, it can be viewed as an injective hommorphism $j: \mathscr{O}_{\mathbb{P}_S^n}(-rd) \to \mathscr{O}_{\mathbb{P}_S^n}$ which remains injective with smooth associated divisor when restricted to any geometric fibre. Then consider the element $(\mathbb{P}_S^n \to S, \mathscr{O}_{\mathbb{P}_S^n}(-rd), j: \mathscr{O}_{\mathbb{P}_S^n}(-rd) \to \mathscr{O}_{\mathbb{P}_S^n})$. The above two functors are pseudo-inverses and hence $\widetilde{\mathscr{D}}_{n,r,d}$ is isomorphic to the scheme $\mathbb{A}_{sm}(n,rd)$. Note that group sheaf $\textbf{\underline{Aut}}(\mathbb{P}_{\mathbb{Z}}^n, \mathscr{O}(-rd))$ which assigns to every scheme $S$ the group of isomorphisms  $\phi: (\mathbb{P}_{\mathbb{Z}}^n, \mathscr{O}(-rd))) \to (\mathbb{P}_{\mathbb{Z}}^n, \mathscr{O}(-rd)))$ is isomorphic to $GL_{n+1}/\mu_{rd}$ (see \cite{AV}, proof of Theorem $4.1$). Moreover since every pair $(P \to S, \mathscr{L})$ where $P \to S$ is Brauer--Severi and $\mathscr{L}$ is a line bundle which restricts to degree $-rd$ in every geometric fibre is \textrm{\'et}ale locally isomorphic to $(\mathbb{P}_S^n, \mathscr{O}(-rd)))$, we have that $\widetilde{\mathscr{D}}_{n,r,d}$ is is a principal bundle over ${\mathscr{D}}_{n,r,d}$. Hence ${\widetilde{\mathscr{D}}_{n,r,d}} \cong \mathbb{A}_{sm}(n,rd)$. Clearly the action is given by $[A] \cdot f(x) = f(A^{-1} x)$ and that concludes the proof.\QEDB
\color{black}

\subsection{Gerbe structure of \texorpdfstring{$\mathscr{H}_{n,r,d}$}{} over \texorpdfstring{$\mathscr{D}_{n,r,d}$}{}}\label{sub2}

We aim to realize $\mathscr{H}_{n,r,d}$ as a gerbe over $\mathscr{D}_{n,r,d}$. 
The following is the main theorem of this section, it shows that the stack of simple cyclic covers over a projective space is a $\mu_r$-gerbe over the stack $\mathscr{D}_{n,r,d}$.

\begin{theorem}\label{stacks of cyclic covers is a gerbe}
There is a morphism $F: \mathscr{H}_{n,r,d} \to \mathscr{D}_{n,r,d}$ of algebraic stacks over an algebraically closed field $\mathbb{k}$, \textcolor{black}{that} realizes $\mathscr{H}_{n,r,d}$ as a gerbe with relative automorphism group sheaf $\mu_r$ 
\end{theorem}

\noindent\textbf{Proof.} We \textcolor{black}{cite} Lemma $8.11.3$ \cite[\href{https://stacks.math.columbia.edu/tag/06NY}{Tag 06NY}]{stacks-project}, ($2$). We note the explicit description of the map $F$. Given an object $ \zeta = (P \to S, \mathscr{L}, i: \mathscr{L}^{\otimes r} \to \mathscr{O}_P)$ of $\mathscr{H}_{n,r,d}$ as in Definition ~\ref{definition of the stack}, $F(\zeta) = (P \to S, \mathscr{L}^{\otimes r}, i: \mathscr{L}^{\otimes r} \to \mathscr{O}_P)$ which is an element of $\mathscr{D}_{n,r,d}$ over $S$. Note that after an \'etale pullback $\{S_i \to S\}$, we have $P|_{S_i} \cong \mathbb{P}_{S_i}^n$. Now condition (a) of Lemma $8.11.3$ follows from the fact that $\mathscr{O}_{\mathbb{P}_S^n}(-rd)$ is $r-$ divisible in $\textrm{Pic}(\mathbb{P}_S^n)$ and condition (b) follows from the fact that an $S$-automorphism $\mathscr{O}_{\mathbb{P}_S^n}(-rd) \to \mathscr{O}_{\mathbb{P}_S^n}(-rd)$ is given by an element of $\mathbb{G}_m(\mathbb{P}_k^n) = k^*$ and $k$ being algebraically closed admits an $r-$th root which induces an $S$-automorphism between
$\mathscr{O}_{\mathbb{P}_S^n}(-d) \to \mathscr{O}_{\mathbb{P}_S^n}(-d)$.

We treat $\mathscr{H}_{n,r,d}$ as a stack over the site $\mathscr{D}_{n,r,d}$ (with topology induced from \'etale topology on the site $(\textrm{Sch/}\mathbb{k})$) and show that the automorphism group of an object $x \to F(x)$ is $\mu_r$. Note that since $\textrm{Aut}(x/F(x))$ consists of automorphisms $\mathscr{L}$ that induces identity on $\mathscr{L}^{\otimes r}$ and $\mathscr{O}_S$, we have that $\textrm{Aut}(x/F(x)) = \mu_r.$\QEDB\par 

\begin{remark}\label{rigidification}
By \cite{MM}, a general smooth polynomial \textcolor{black}{in $n+1$ variables} of degree $k \geq 4$ has \textcolor{black}{stabilizer} group $\mu_k$ \textcolor{black}{under the action of $GL_{n+1}$ given by $A \cdot f(x) = f(A^{-1}(x))$}. Hence the generic automorphism group of $\mathscr{H}_{n,r,d}$ is $\mu_r$ and by Theorem \ref{orbifold}, $\mathscr{D}_{n,r,d}$ has generic trivial stabilizer. Since by Theorem ~\ref{stacks of cyclic covers is a gerbe} $\mathscr{H}_{n,r,d}$ is an abelian gerbe over $\mathscr{D}_{n,r,d}$, the latter is the $\mu_r$-rigidification of $\mathscr{H}_{n,r,d},$ \textcolor{black}{ where for the definition of rigidification see \cite{AOV}, Appendix C}. Since smooth hypersurfaces of degree at least four are stable under the action of $GL_{n+1}/\mu_{rd}$ or $GL_{n+1}/\mu_d$, we have that the quotient stacks have finite stabilizers and hence by \cite{E}, Theorem $4.7$ and Theorem $4.18$, we have that  both are separated DM stacks (recall that a DM stack is {\bf separated} if and only if the diagonal is finite by Remark $4.16$ of \cite{E}). The stability of the hypersurfaces also show that there exists a geometric quotient (see \cite{Mum}) of $\mathbb{A}_{sm}(n,rd)$ by either $GL_{n+1}/\mu_{d}$ or $GL_{n+1}/\mu_{rd}$. Since $\mathscr{D}_{n,r,d}$ is a $\mu_r$-rigidification of $\mathscr{H}_{n,r,d}$ they both have the same coarse moduli space. But by \cite{E}, Proposition $4.26$ the respective coarse moduli are the respective geometric quotients and hence the geometric quotients are isomorphic 
\begin{equation*}
    \mathbb{A}_{sm}(n,rd)//(GL_{n+1}/\mu_d) \cong \mathbb{A}_{sm}(n,rd)//(GL_{n+1}/\mu_{rd}).
\end{equation*} 
We denote by $M_{n,r,d}$ the common coarse moduli of the stacks. Since $\mathscr{D}_{n,r,d}$ has generic trivial stabilizer, there exist an open substack $\mathscr{D}_{n,r,d}^0 \subset \mathscr{D}_{n,r,d}$ and an open subscheme $M_{n,r,d}^0 \subset M_{n,r,d}$ such that $\mathscr{D}_{n,r,d}^0 \cong M_{n,r,d}^0$. 
\end{remark}

Theorem \ref{stacks of cyclic covers is a gerbe} immediately shows the following

\begin{corollary}\label{grb0}
Let $\mathscr{H}_{n,r,d}^0$ be the open substack of $\mathscr{H}_{n,r,d}$ which is the inverse image of $M_{n,r,d}^0$ under the coarse moduli map $p: \mathscr{H}_{n,r,d} \to M_{n,r,d}$. Then $p: \mathscr{H}_{n,r,d}^0 \to M_{n,r,d}^0$ realizes $\mathscr{H}_{n,r,d}^0$ as a $\mu_r$-gerbe over the scheme $M_{n,r,d}^0$. 
\end{corollary}

\subsection{Degree of line bundles on the universal Brauer--Severi scheme \texorpdfstring{$(P \to M_{n,r,d}^0)$}{}}\label{sub3}

In this section, we aim to study the least possible fibrewise degree of a line bundle on the universal Brauer--Severi scheme $(P \to M_{n,r,d}^0)$. \textcolor{black}{Note that for a cyclic group $\mu_k$, there is a natural map from $GL_{n+1}/\mu_k \to PGL_{n+1}$ whose fibres are $\mathbb{G}_m/\mu_k$. The natural action of $\mathbb{G}_m/\mu_k$ on $GL_{n+1}/\mu_k$ which is multiplication by equivalence classes of scalar matrices preserves the map $GL_{n+1}/\mu_k \to PGL_{n+1}$ and acts freely and transitively on the fibres. Also the pullback of $GL_{n+1}/\mu_k \to PGL_{n+1}$ by the same map $GL_{n+1}/\mu_k \to PGL_{n+1}$ yields $GL_{n+1}/\mu_k \times \mathbb{G}_m/\mu_k \to GL_{n+1}/\mu_k$. Hence $GL_{n+1}/\mu_k \to PGL_{n+1}$ is a map of $\mathbb{G}_m/\mu_k$ torsors or principal $\mathbb{G}_m/\mu_k$ bundles. Since we have an isomorphism from $\mathbb{G}_m/\mu_k \cong \mathbb{G}_m$ induced by the map $\mathbb{G}_m \to \mathbb{G}_m$ given by $\lambda \to \lambda^k$, the natural map from $GL_{n+1}/\mu_k \to PGL_{n+1}$ realizes $GL_{n+1}/\mu_k$ as a $\mathbb{G}_m$ torsor over $PGL_{n+1}$.} 

\begin{lemma}\label{torsor morphisms}
There exists a map $GL_{n+1}/\mu_d \to GL_{n+1}/\mu_{d'}$ of $\mathbb{G}_m$ torsors over $PGL_{n+1}$ if and only if there exists an integer $k$ such that $d \mid k(n+1)+d'$, i.e., $\textrm{gcd}(d,n+1) \mid d'$.
\end{lemma}

\noindent\textbf{Proof.} Suppose $d \mid k(n+1)+d'$. Then consider the map $\varphi:$ $GL_{n+1}/\mu_d \to GL_{n+1}/\mu_{d'}$ of $\mathbb{G}_m$ torsors over $PGL_{n+1}$ given by 
\begin{equation*}
    \varphi([A]) = [\textrm{det}(A)^{\frac{k}{d'}}A].
\end{equation*}
\textcolor{black}{for some choice of $d'$-th root of $\textrm{det}(A)^k$. The class of $[\textrm{det}(A)^{\frac{k}{d'}}A]$ is well-defined in  $GL_{n+1}/\mu_{d'}$.}   We show that $\varphi$ is well defined. Suppose that $\zeta_d \in \mu_d$. Then 
 \begin{equation*}
    \varphi([\zeta_d]) = [\textrm{det}([\zeta_d])^{\frac{k}{d'}}[\zeta_d]] 
     = \left[\zeta_d^{\frac{k(n+1)}{d'}}[\zeta_d]\right] 
     = \left[\zeta_d^{\frac{k(n+1)+d'}{d'}}\right] 
     = \left[\zeta_d^{\frac{ad}{d'}}\right].
\end{equation*}
\textcolor{black}{where $ad = k(n+1)+d'$}. Now, since $\zeta_d^{\frac{ad}{d'}}$ is a $d'$-th root of unity we have that $\varphi$ is well-defined. The fact that $\varphi$ is a homomorphism can easily be seen from the following
\begin{equation*}
    \varphi([AB])  = [\textrm{det}(AB)^{\frac{k}{d'}}AB] 
     = [\textrm{det}(A)^{\frac{k}{d'}}\textrm{det}(B)^{\frac{k}{d'}}AB] 
     = [\textrm{det}(A)^{\frac{k}{d'}}A][\textrm{det}(B)^{\frac{k}{d'}}B] 
     = \varphi([A])\varphi([B]).
\end{equation*}

Conversely, suppose that $\varphi: GL_{n+1}/\mu_d \to GL_{n+1}/\mu_{d'}$ be a morphism of $\mathbb{G}_m$ torsors over $PGL_{n+1}$. Then $\varphi([A]) = [f(A)A]$ where $f(A) \in \mathbb{G}_m$. Note that, since $\varphi$ is a homomorphism of groups, we have that $$\varphi([AB]) = [f(AB)AB] = \varphi([A])\varphi([B]) = [f(A)A][f(B)B] = [f(A)f(B)AB].$$
Hence $f(AB) = f(A)f(B) \zeta_{d'}$ where $\zeta_{d'}$ is a $d'$-th root of unity and therefore $\overline{f}: GL_{n+1} \to \mathbb{G}_m/\mu_{d'}$ defined by $\overline{f}(A) = [f(A)]$ is a well-defined homomorphism. Now consider the isomorphism $\psi: \mathbb{G}_m/\mu_{d'} \to \mathbb{G}_m$ given by $x \to x^{d'}$. Composing with $\psi$ we have a homomorphism $g = \overline{f}^{d'}: GL_{n+1}/\mu_d \to \mathbb{G}_m$. Hence $[f(A)]^{d'} = \textrm{det}(A)^k$ for some $k\in\mathbb{Z}$ and hence $\overline{f}(A) = [\textrm{det}(A)^{\frac{k}{d'}}]$. Hence $$\varphi([A]) = [f(A)A] = [f(A)][A] = \overline{f}(A)[A] = [\textrm{det}(A)^{\frac{k}{d'}}A].$$ We now find out the condition imposed by $\varphi$ being well-defined. Let $\zeta_d$ is a $d$-th root of unity and notice 
\begin{equation*}
    \varphi([\zeta_d]) = [\textrm{det}([\zeta_d])^{\frac{k}{d'}}[\zeta_d]] 
     = \left[\zeta_d^{\frac{k(n+1)}{d'}}[\zeta_d]\right] 
     = \left[\zeta_d^{\frac{k(n+1)+d'}{d'}}\right]. 
\end{equation*}
Thus, since $\varphi$ is well-defined, we must have that $\zeta_d^{k(n+1)+d'} = 1$ which implies $d | k(n+1)+d'$. \QEDB

\vspace{5pt}

The 
following lemma (see also \cite{GV1}))
is essential to this article.

\begin{lemma}\label{relative picard of universal Brauer Severi}
Let the unique section $p$ to $\mathscr{D}_{n,r,d}^0 \to M_{n,r,d}^0$ be denoted by the object $(P \to M_{n,r,d}^0, D)$ where $P \to S$ is Brauer--Severi and $D$ is a divisor of degree $rd$ on every geometric fibre. Let $U \subseteq M_{n,r,d}^0$ be open. Then $$\textrm{Pic} (P|_U/U) \cong \textrm{Pic}^{GL_{n+1}/\mu_{rd}}(\mathbb{P}^n)/K.$$
\textcolor{black}{where $K \subseteq \textrm{Lin}^{GL_{n+1}/\mu_{rd}}(\mathscr{O}_{\mathbb{P}^n})$ and $\textrm{Lin}^{GL_{n+1}/\mu_{rd}}(\mathscr{O}_{\mathbb{P}^n}) = \hat{G}$ is the subgroup of pairs $(\mathscr{O}_{\mathbb{P}^n}, f) $ where $f$ is a $GL_{n+1}/\mu_{rd}$ linearization of $\mathscr{O}_{\mathbb{P}^n}$.}
\end{lemma}

\noindent\textbf{Proof.} We write the proof taking $U = M_{n,r,d}^0$ and note that same proof goes for any $U \subseteq M_{n,r,d}^0$. Let $\mathbb{A}_{sm}^0(n,rd)\subset \mathbb{A}_{sm}(n,rd)$ be the open set consisting of smooth forms with trivial stabilizer inside $GL_{n+1}/\mu_{rd}$. Further, let $\mathbb{D}$ be the incidence divisor inside $\mathbb{A}_{sm}^0(n,rd) \times \mathbb{P}^n$ defined by $$\mathbb{D}:= \left\{(x,y) \in \mathbb{A}_{sm}^0(n,rd) \times \mathbb{P}^n \mid y \in V_{+}(x) \right\}.$$ There is a natural free action of $GL_{n+1}/\mu_{rd}$ on $\mathbb{A}^0_{sm}(n,rd)$ such that $M_{n,r,d}^0$ is its geometric quotient. Let $GL_{n+1}/\mu_{rd}$ act diagonally on $\mathbb{A}_{sm}^0(n,rd) \times \mathbb{P}^n$ so that the action is free and admits a geometric quotient $P'$. Now clearly $\mathbb{D}$ is invariant under the action of $GL_{n+1}/\mu_{rd}$ and the induced action of $GL_{n+1}/\mu_{rd}$ is free. Hence there  exists a geometric quotient $D'$. Consider the commutative diagram with the left vertical chain consisting $GL_{n+1}/\mu_{rd}$ equivariant maps and the horizontal maps being geometric quotients under the free action of $GL_{n+1}/\mu_{rd}$.  
\[
\begin{tikzcd}
\mathbb{D}  \arrow{r} & \mathbb{A}_{sm}^0(n,rd) \times \mathbb{P}^n \arrow[r] \arrow[d] & P' \arrow[d] \\ 
& \mathbb{A}_{sm}^0(n,rd) \arrow[r] & M_{n,r,d}^0 
\end{tikzcd}
\]
Now the incidence divisor $\mathbb{D}$ is in the linear system of the $GL_{n+1}/\mu_{rd}$ linearized line bundle $\mathscr{O}_{\mathbb{P}_{\mathbb{A}_{sm}^0(n,rd)}^n}(rd)$. Since the action of $GL_{n+1}/\mu_{rd}$ on $\mathbb{A}_{sm}^0(n,rd) \times \mathbb{P}^n$ is free, $P'$ is actually the quotient stack $$[\mathbb{A}_{sm}^0(n,rd) \times \mathbb{P}^n/GL_{n+1}/\mu_{rd}]$$ and hence $D'$ determines a divisor on $P'$ of relative degree $rd$ on each fibre over a geometric point of $M_{n,r,d}^0$. Hence we have the following diagram where horizontal rows are geometric quotients of free actions of $GL_{n+1}/\mu_{rd}$.
\[
\begin{tikzcd}
  \mathbb{D}  \arrow[r]\arrow{d} & D  \arrow[d]    \\
 \mathbb{A}_{sm}^0(n,rd) \times \mathbb{P}^n \arrow[r] \arrow[d] & P' \arrow[d] \\ 
\mathbb{A}_{sm}^0(n,rd) \arrow[r] & M_{n,r,d}^0 
\end{tikzcd}
\]
Hence the object $(P' \to M_{n,r,d}^0, D)$ induces a section $M_{n,r,d}^0 \to \mathscr{D}_{n,r,d}^0$. But by uniqueness of the section, we have that $(P' \to M_{n,r,d}^0, D') = (P \to M_{n,r,d}^0, D)$. Hence 
\begin{align*}
    \textrm{Pic}(P) = \textrm{Pic}(P') =  \textrm{Pic}^{GL_{n+1}/\mu_{rd}}(\mathbb{A}_{sm}^0(n,rd) \times \mathbb{P}^n) 
\end{align*}
\color{black}
Thus, $\textrm{Pic}(P/M_{n,r,d}^0) = \textrm{Pic}^{GL_{n+1}/\mu_{rd}}(\mathbb{A}_{sm}^0(n,rd) \times \mathbb{P}^n)/\pi^*(\textrm{Pic}^{GL_{n+1}/\mu_{rd}}(\mathbb{A}_{sm}^0(n,rd)))$ where $\pi: \mathbb{A}_{sm}^0(n,rd) \times \mathbb{P}^n) \to \mathbb{A}_{sm}^0(n,rd)$ is the projection under the first coordinate. For the rest of the proof, for simplicity of notation let $G = GL_{n+1}/\mu_{rd}$, $X = \mathbb{A}_{sm}^0(n,rd)$ and $Y = \mathbb{P}^n$. Since atleast one of $X$ and $Y$ is rational (in fact both) and both being normal we have that $\textrm{Pic}(X \times Y) = \textrm{Pic}(X) \times \textrm{Pic}(Y)$. Also note that $\textrm{Pic}(X) = 0$. Now consider the following exact sequence (see \cite{D03}, Theorem $7.2$) 
\begin{equation*}
    0 \to \textrm{Lin}^G(\mathscr{O}_{X \times Y}) \to \textrm{Pic}^{G}(X \times Y) \to \textrm{Pic}(X \times Y) = \textrm{Pic}(X) \times \textrm{Pic}(Y) = \textrm{Pic}(Y) 
\end{equation*}
where $\textrm{Lin}^G(\mathscr{O}_{X \times Y})$ is the subgroup of $\textrm{Pic}^{G}(X \times Y)$ consisting of pairs $(\mathscr{O}_{X \times Y},f)$ where $f$ is a $G$ linearization of $\mathscr{O}_{X \times Y}$.
Similarly there is an exact sequence 
\begin{equation*}
    0 \to \textrm{Lin}^G(\mathscr{O}_{X}) \times \textrm{Lin}^G(\mathscr{O}_{Y}) \to \textrm{Pic}^{G}(X) \times \textrm{Pic}^{G}(Y) \to  \textrm{Pic}(X) \times \textrm{Pic}(Y) = \textrm{Pic}(Y)  
\end{equation*}
Now note that since $\textrm{Pic}(X) = 0$, we have that the image of the map $\textrm{Pic}^{G}(X) \times \textrm{Pic}^{G}(Y) \to  \textrm{Pic}(Y)$ is $H \subset \textrm{Pic}(Y)$ where $H$ is the image of the map $\textrm{Pic}^{G}(Y) \to \textrm{Pic}(Y)$. Now using the facts that $G$ acts diagonally on $X \times Y$ and $\textrm{Pic}(X) = 0$, we claim that the image of the map $\textrm{Pic}^{G}(X \times Y) \to \textrm{Pic}(Y)$ is also $H$. First of all, it is clear that $H$ is in the image of $\textrm{Pic}^{G}(X \times Y) \to \textrm{Pic}(Y)$. Consider the commutative diagrams

\vspace{5pt}

\begin{minipage}{.5\textwidth}
\centering
\begin{tikzcd}
   G \times X \times Y \arrow[r, "\pi_2"] \arrow{d}{\sigma}   & G \times Y \arrow{d}{\sigma_2}  \\
  X \times Y \arrow[r, "p_2"] & Y 
   \end{tikzcd}
\end{minipage}
\begin{minipage}{.4\textwidth}
\centering
   \begin{tikzcd}
  G \times X \times Y \arrow[r, "\pi_2"] \arrow{d}{\pi}   & G \times Y \arrow{d}{\lambda}  \\
   X \times Y \arrow[r, "p_2"] & Y 
  \end{tikzcd} 
\end{minipage}

\noindent where $\sigma$ and $\sigma_2$ are the group actions on $X \times Y$ and $Y$ respectively and $\pi$ and $\lambda$ are projections onto the respective coordinates. Suppose $L \in \textrm{Pic}(Y)$ is in the image of the map $\textrm{Pic}^{G}(X \times Y) \to \textrm{Pic}(Y)$. Since $\textrm{Pic}(X) = 0$, we have an isomorphism $$ \pi^*(p_2^*(L))) \cong \sigma^*(p_2^*(L))).$$
But $\pi^*(p_2^*(L))) = \pi_2^*(\lambda^*(L)))$ and $\sigma^*(p_2^*(L))) \cong \pi_2^*(\sigma_2^*(L)))$ using the fact the action is diagonal. Hence we have an isomorphism $$ \pi_2^*(\lambda^*(L)) \cong \pi_2^*(\sigma_2^*L).$$
Now considering that $\textrm{Pic}(G \times X \times Y) = \textrm{Pic}(G \times Y) \times \textrm{Pic}(X) = \textrm{Pic}(G \times Y)$, $\pi_2^*$ is injective (and in this case it is an isomorphism), and we have that $$\lambda^*(L) \cong \sigma_2^*(L)$$ which imples that $L$ has a $G$ linearization. Hence our claim is proved. 

So we have a commutative diagram of exact sequences
\[
\begin{tikzcd}
  0  \arrow[r] & \textrm{Lin}^G(\mathscr{O}_{X}) \times \textrm{Lin}^G(\mathscr{O}_{Y}) \arrow[r] \arrow{d}{a}  & \textrm{Pic}^{G}(X) \times \textrm{Pic}^{G}(Y) \arrow[r]\arrow{d}{b} & H \arrow[r] \arrow[d, "\sim"] & 0   \\
 0  \arrow[r] & \textrm{Lin}^G(\mathscr{O}_{X \times Y}) \arrow[r] & \textrm{Pic}^{G}(X \times Y) \arrow[r] & H \arrow[r] & 0
\end{tikzcd}
\]
We first show that $a$ is surjective. Consider the section $\textrm{Lin}^G(\mathscr{O}_{X}) \to \textrm{Lin}^G(\mathscr{O}_{X}) \times \textrm{Lin}^G(\mathscr{O}_{Y})$ given by the trivial $G$ linearization $\mathscr{O}_Y$. It is enough to show that the map $$\textrm{Lin}^G(\mathscr{O}_{X}) \to \textrm{Lin}^G(\mathscr{O}_{X \times Y})$$ obtained by composing $$\textrm{Lin}^G(\mathscr{O}_{X}) \to \textrm{Lin}^G(\mathscr{O}_{X}) \times \textrm{Lin}^G(\mathscr{O}_{Y}) \to \textrm{Lin}^G(\mathscr{O}_{X \times Y})$$ where the first map is the section induced by the trivial linearization on $\mathcal{O}_Y$, is surjective. To see this note that by Rosenlicht's theorem, there is a surjection from $\mathscr{O}(G)^* \otimes \mathscr{O}(X)^* \to \mathscr{O}(G \times X)^*$. Similarly there is a surjection from $\mathscr{O}(X)^* \otimes \mathscr{O}(G \times Y)^* \to \mathscr{O}(G \times X \times Y)^*$. But since $Y$ is proper and connected, we have that $\mathscr{O}(G \times Y)^* \cong \mathscr{O}(G)^*$. So we have a surjection from
$\mathscr{O}(X)^* \otimes \mathscr{O}(G)^* \to \mathscr{O}(G \times X \times Y)^*$. But the last map factors through $\mathscr{O}(G \times X)^*$ and hence we have a surjection from 
$\mathscr{O}(G \times X)^* \to \mathscr{O}(G \times X \times Y)^*$. We need to show that this induces a surjection from $\textrm{Lin}^G(\mathscr{O}_{X}) \to \textrm{Lin}^G(\mathscr{O}_{X \times Y})$. We once again use the fact $G$ acts diagonally on $X \times Y$. Suppose an element $\psi \in \mathscr{O}(G \times X \times Y)^*$ satisfies the cocycle condition $$ \psi(gg',x,y) = \psi(g, g'x, g'y)\psi(g'x,g'y).$$ Consider the preimage $\phi \in \mathscr{O}(G \times X)^*$ which implies $\psi = \textrm{pr}_1^*(\phi)$ where $\textrm{pr}_1: G \times X \times Y \to G \times X$ is the projection onto the first two coordinates. Then 
$$\textrm{pr}_1^*(\phi)(gg',x,y) = \textrm{pr}_1^*(\phi)(g, g'x, g'y)\textrm{pr}_1^*(\phi)(g'x,g'y),$$
which implies
$$\phi (\textrm{pr}_1(gg',x,y)) = \phi (\textrm{pr}_1(g, g'x, g'y))\phi (\textrm{pr}_1(g'x,g'y)),$$
which implies
$$\phi (gg',x) = \phi (g, g'x)\phi (g',x).$$
Hence $\phi$ also satisfies the cocycle condition. So we have that $a$ is surjective. 

Now apply snake lemma to the above commutative diagram of exact sequences to obtain
$$0 \to \textrm{Ker}(a) \to \textrm{Ker}(b) \to 0 \to 0 \to \textrm{Coker}(b) \to 0.$$
Therefor $\textrm{Ker}(b) \cong \textrm{Ker}(a) $ and $b$ is surjective. Hence we have an exact sequence
$$ 0 \to \textrm{Ker}(a) \to \textrm{Pic}^{G}(X) \times \textrm{Pic}^{G}(Y) \to \textrm{Pic}^{G}(X \times Y) \to 0.  $$
Now consider the following commutative diagram of exact sequences
\[
\begin{tikzcd}
  0  \arrow[r] &  0 \arrow[r] \arrow{d}  & \textrm{Pic}^{G}(X) \arrow[r]\arrow{d}{c} & \textrm{Pic}^{G}(X) \arrow[r] \arrow[d, "d"] & 0   \\
 0  \arrow[r] & \textrm{Ker}(a) \arrow[r] & \textrm{Pic}^{G}(X) \times \textrm{Pic}^{G}(Y) \arrow[r] & \textrm{Pic}^{G}(X \times Y) \arrow[r] & 0
\end{tikzcd}
\]
where the map $c$ is the section induced by the identity element of $\textrm{Pic}^G(Y)$ and $d$ is the map induced by the projection $X \times Y \to X$. Then by snake lemma we have
$$ 0 \to \textrm{Ker}(c) \to \textrm{Ker}(d) \to \textrm{Ker}(a) \to \textrm{Pic}^G(Y) \to \textrm{Coker}(d) = \textrm{Pic}(P/M_{n,r,d}^0) \to 0. $$
But by the naturality of all maps involved, the map from $\textrm{Ker}(a) \to \textrm{Pic}^G(Y)$ factors as $$\textrm{Ker}(a) \hookrightarrow \textrm{Lin}^G(\mathscr{O}_{X}) \times \textrm{Lin}^G(\mathscr{O}_{Y}) \to \textrm{Lin}^G(\mathscr{O}_{Y}) \hookrightarrow \textrm{Pic}^G(Y).$$ Thus the image of the map $\textrm{Ker}(a) \to \textrm{Pic}^G(Y)$ lands in $\textrm{Lin}^G(\mathscr{O}_{Y})$ and hence the lemma is proved.\QEDB

\vspace{5pt}

The following lemma computes the least possible degree of a line bundle on the Brauer--Severi scheme associated to the universal rational section $M_{n,r,d}^0 \to \mathscr{D}_{n,r,d}^0$. This in particular determines whether the universal Brauer--Severi scheme $P \to M_{n,r,d}^0$ is a projective bundle or not.

\begin{lemma}{\label{least linearized degree}}
Let \[
s:= \min\Set{ \vert t\vert | \begin{array}{l} \exists \text{ a line bundle $\mathscr{L}$ on the universal Brauer--Severi scheme}\\ \text{$P \to M_{n,r,d}^0$ of degree $t \neq 0$ on every geometric fibre}\end{array}}.
\] Then $s = \textrm{gcd}(rd, n+1)$.
\end{lemma}
\noindent\textbf{Proof.} By Lemma ~\ref{relative picard of universal Brauer Severi}, 
\color{black}
\begin{equation*}
    s= \min\Set{\vert t\vert | \exists\text{ a line bundle $\mathscr{L}$ on $\mathbb{P}^n$ of degree $t \neq 0$ such that $\mathscr{L} \in \textrm{Pic}^{GL_{n+1}/\mu_{rd}}(\mathbb{P}^n)/K$}}.
\end{equation*}
where $K \subseteq \textrm{Lin}^{GL_{n+1}/\mu_{rd}}(\mathscr{O}_{\mathbb{P}^n})$. Since $K \subseteq \textrm{Lin}^{GL_{n+1}/\mu_{rd}}(\mathscr{O}_{\mathbb{P}^n})$ we have that 
\begin{equation*}
    s= \min\Set{\vert t\vert | \exists\text{ a line bundle $\mathscr{L}$ on $\mathbb{P}^n$ of degree $t \neq 0$ such that $\mathscr{L}$ admits a $GL_{n+1}/\mu_{rd}$ linearization}}.
\end{equation*}
\color{black}
 Since a $GL_{n+1}/\mu_{rd}$ linearization of $\mathscr{O}(t)$ implies a $GL_{n+1}/\mu_{rd}$ linearization of $\mathscr{O}(-t)$, we can assume $t > 0$. Since Aut ($\mathbb{P}^n, \mathscr{O}(t)) = GL_{n+1}/\mu_t$, we have that by definition $s \leq rd$. Now if $t = n+1$, then by Lemma ~\ref{torsor morphisms}, we have that there is a morphism of $\mathbb{G}_m$ torsors over $PGL_{n+1}$ from $GL_{n+1}/\mu_{rd}$ to $GL_{n+1}/\mu_{n+1}$ (since $\textrm{gcd}(rd,n+1) \mid n+1$) and hence there exists a $GL_{n+1}/\mu_{rd}$ linearization of $\mathscr{O}(n+1)$. Therefore there exists a $GL_{n+1}/\mu_{rd}$ linearization of $\mathscr{O}(\textrm{gcd}(rd, n+1))$. Hence $s \leq \textrm{gcd}(rd, n+1)$. Now suppose that there exists a $GL_{n+1}/\mu_{rd}$ linearization of $\mathscr{O}(t)$. Then there exists a map of $\mathbb{G}_m$ torsors over $PGL_{n+1}$ from $GL_{n+1}/\mu_{rd} \to GL_{n+1}/\mu_{t}$. Hence $rd \mid k(n+1) + t$ for some integer $k$. Then $\textrm{gcd}(rd, n+1) \mid t$. Since both are positive we conclude $$t \geq \textrm{gcd}(rd, n+1).$$ Hence $s = \textrm{gcd}(rd, n+1)$.\QEDB

\subsection{
Tautological family on Zariski open sets of the coarse moduli space of \texorpdfstring{$\mathscr{H}_{n,r,d}$}{}}\label{sub4} This section is devoted to the proof of Theorem \ref{main1intro}.

\subsubsection{Existence of tautological families}

The following is the main theorem on the existence of tautological sections over an open subset of $M_{n,r,d}$. This together with Lemma $99.27.6$ (\cite[\href{https://stacks.math.columbia.edu/tag/06QB}{Tag 06QB}]{stacks-project}) determines whether the gerbe (see Corollary \ref{grb0}) $\mathscr{H}_{n,r,d}^0 \to M_{n,r,d}^0$ is trivial over a Zariski open set of $M_{n,r,d}$ or not.

\begin{theorem}\label{tautological family} 
The following statements hold true.
\begin{itemize}
    \item[(1)] There exists a tautological family of cyclic covers over an open set $U$ of the coarse moduli $M_{n,r,d}$ if and only if $\textrm{gcd}(rd, n+1) \mid d$.
    \item[(2)] When such a family exists, the Brauer--Severi scheme $P \to U$ associated to the tautological family is trivial if and only if $\textrm{gcd}(rd, n+1) = 1$.
    \item[(3)] \textcolor{black}{Assume that the codimension of the complement of $M_{n,r,d}^0$ in $M_{n,r,d}$ is greater than or equal to $2$}.  The family when it exists is unique when restricted to the complement of all hyperplane sections in $M_{n,r,d}^0$ in the linear systems of each unique line bundle of order $k$ where $k \mid \textrm{gcd}(r,n+1)$. In particular, if $\textrm{gcd}(r,n+1) = 1$, then for those open sets over which a tautological family exists, such a family is unique.
\end{itemize}
\end{theorem}

\color{black}
\noindent\textbf{Proof.} (1) Let the unique section to the natural map $p: \mathscr{D}_{n,r,d}^0 \to M_{n,r,d}^0$ be denoted by $(\pi: P \to M_{n,r,d}^0, D)$ where $P \to M_{n,r,d}^0$ is a Brauer--Severi scheme of relative dimension $n$ and $D$ is a divisor smooth over $M_{n,r,d}^0$ with relative degree $rd$. 

Now suppose that $\textrm{gcd}(rd, n+1) \mid d$. By Lemma ~\ref{torsor morphisms} there exist a morphism from $GL_{n+1}/\mu_{rd} \to GL_{n+1}/\mu_d$. Hence $\mathscr{O}_{\mathbb{P}^n}(d)$ has a $GL_{n+1}/\mu_{rd}$ linearization and consequently $\mathscr{O}_{\mathbb{P}^n}(d) \in \textrm{Pic}^{GL_{n+1}/\mu_{rd}}(\mathbb{P}^n)$. \textcolor{black}{Consequently, the class $[\mathscr{O}_{\mathbb{P}^n}(d)] \in \textrm{Pic}(P/ M_{n,r,d}^0)$ is non-zero}  
we have a line bundle $\mathscr{L}$ on $P$ with \textcolor{black}{vertical} degree $d$. Then the divisor $D$ belongs to the linear system of the line bundle $\mathscr{L}^{\otimes r} \otimes \pi^*(M)$ for some line bundle $M$ on $M_{n,r,d}^0$. Then, after possibly shrinking $M_{n,r,d}^0$, we can assume that $D$ is in the linear system of an $r$ divisible line bundle on $P$. So one can construct an $r$ cyclic cover over $P \to S$ branched along $D$ to get a tautological family of cyclic covers over $\mathbb{P}^n$. 

Now suppose that there exists a tautological family of cyclic covers over an open set $U$ of $M_{n,r,d}$. Intersecting with $M_{n,r,d}^0$ we can assume that $U$ is contained inside $M_{n,r,d}^0$. This implies that the divisor $D|_U$ on $\pi: P|_U \to U$ is contained in the linear system of an $r$ divisible line bundle $\mathscr{L}^{\otimes r}$. Now by Lemma ~\ref{relative picard of universal Brauer Severi} we have that $\mathscr{L} \in \textrm{Pic}(P|_U/U) = \textrm{Pic}^{GL_{n+1}/\mu_{rd}}(\mathbb{P}^n)/K$. \textcolor{black}{Since $\mathscr{L}$ has vertical degree $d$, we conclude that} $\mathscr{O}(d)$ has a linearization under the action of $GL_{n+1}/\mu_{rd}$ \textcolor{black}{and hence} we have that there exists a homomorphism of $\mathbb{G}_m$ torsors over $PGL_{n+1}$ between $GL_{n+1}/\mu_{rd} \to GL_{n+1}/\mu_{d}$ and our result follows from Lemma ~\ref{torsor morphisms}. \par 
 
(2). Note that a Brauer--Severi scheme $P$ is Zariski locally trivial if and only if there exists a line bundle of degree one when restricted to a fibre. Now our second assertion follows from Lemma ~\ref{least linearized degree}. \par

(3) By \cite{GV1}, Lemma $3.2$, part $\textrm{(ii)}$, one has that if a family exists over an open set $U \subset M_{n,r,d}^0$, the number of such families is the number of $r$ torsion elements of $\textrm{Pic}(U)$. Now note that by excision there is a surjection  $\textrm{Pic}(M_{n,r,d}^0) \to \textrm{Pic}(U)$
where the kernel consists of all hypersurfaces in $M_{n,r,d}^0-U$. \textcolor{black}{Since by our assumption, the codimension of the complement of $M_{n,r,d}^0$ in $M_{n,r,d}$ is greater than or equal to 2, $\textrm{Pic}(M_{n,r,d}^0) = \textrm{Pic}(\mathscr{D}_{n,r,d}^0) = \textrm{Pic}(\mathscr{D}_{n,r,d})$} 
\textcolor{black}{and hence} \textrm{Pic}($M_{n,r,d}^0$) is cyclic of order $(rd-1)^n
\textrm{gcd}(rd, n+1)$ \textcolor{black}{by Lemma ~\ref{comparison of picard groups}, $(3)$.} So the $r-$ torsion elements correspond to hyperplane sections of the unique line bundle in $\textrm{Pic}(M_{n,r,d}^0)$ of order $k$ dividing $\textrm{gcd}(r, (rd-1)^n \textrm{gcd}(rd, n+1)) = \textrm{gcd}(r, n+1)) $. \QEDB

\vspace{5pt}

We now recall the definition of rationality of stack.

\begin{definition}\label{rationality of a stack}(see \cite{BH08}, Section $4$)
\begin{itemize}
    \item[(i)] A stack $\mathscr{X}$ over an algebraically closed field $\mathbb{k}$ is called {\bf rational} if it is birational ($1$- isomorphism between non-empty open substacks) to $\mathbb{A}^n \times BG$ for some group scheme $G$ over $\mathbb{k}$.
    \item[(ii)] A stack $\mathscr{X}$ over an algebraically closed field $\mathbb{k}$ is called {\bf unirational} if it admits a dominant $1$-morphism from a dense open subscheme of $\mathbb{A}^n$ for some $n$.

\end{itemize} 
\end{definition}


\begin{remark}
    It follows from Example ~\ref{toricdm} that a finite abelian gerbe $\mathscr{X}$ is rational if and only it is birational to a Toric Deligne-Mumford stack.
\end{remark}

\begin{example}\label{toricdm}
\textbf{Toric Deligne-Mumford stacks:} a \textit{Deligne-Mumford Torus} $\mathscr{T}$ is isomorphic as a \textit{Picard stack} to 
$T\times BG$ where $T$ is a torus and $G$ is a finite abelian group (\cite{FMN} Proposition $2.6$).
A smooth \textit{toric Deligne-Mumford stack} is a smooth separated Deligne-Mumford stack $X$ together with an open immersion of a Deligne-Mumford torus $\mathscr{T}\xrightarrow[]{i}X$ with dense image such that the action of $\mathscr{T}$ on itself extends to an action $a:\mathscr{T}\times X\rightarrow X.$ (\cite{FMN}, Definition $3.1$). Hence it is clear that a smooth Toric Deligne-Mumford stack $X$ is a rational stack. An important class of Toric Deligne-Mumford stacks given by the Weighted Projective Stack (\cite{FMN} Example 7.27) that we describe below.\\
\noindent\underline{\textit{Weighted Projective Stack.}} Let $$w=(w_{1},w_{2}..,w_{n})\in \mathbb{N}^{n}_{>0},\quad \textrm{gcd}(w_{1},w_{2},,w_{n})=d.$$
The weighted projective stack is given by $$\mathbb{P}(w)=[\mathbb{k}^{n}-(0)/\mathbb{k}^{*}],$$ where the action of $\mathbb{k}^{*}$ on $\mathbb{k}^{n}-(0)$ is as follows
$$(\lambda,(x_{1},x_{2},..,x_{n}))=(\lambda^{w_{1}}x_{1},\lambda^{w_{2}}x_{2},..,\lambda^{w_{n}}x_{n}).$$ Here we denote the $(x_{1},x_{2},..,x_{n})$ as choice of coordinates of $\mathbb{k}^{n}$.
\end{example}

Denote by $C(n,k)$ the coarse moduli parametrizing isomorphism classes of $\textrm{deg}(k)$ hypersurfaces in $\mathbb{P}^n$. The following lemma relates the coarse moduli of $\mathscr{H}_{n,r,d}$ and $\mathscr{D}_{n,r,d}$ to the moduli of hypersurfaces in projective space.

\begin{lemma}\label{Coarse moduli of plane curves 1}
The coarse moduli space of $\mathscr{H}_{n,r,d}$ is isomorphic to the coarse moduli of degree $rd$ hypersurfaces in $\mathbb{P}^n$, i.e., $M_{n,r,d} = \mathbb{P}(n+1,rd) // PGL_{n+1} = C(n,rd)$ .  
\end{lemma}

\noindent\textbf{Proof.} By Remark ~\ref{rigidification}, the coarse moduli of $\mathscr{H}_{n,r,d}$ is isomorphic to the coarse moduli of $\mathscr{D}_{n,r,d}$. Hence it is enough to show that $\mathbb{A}_{sm}(n,rd)//(GL_{n+1}/\mu_{rd}) \cong \mathbb{P}(n+1,rd) // PGL_{n+1}$. The natural map $f: \mathbb{A}(n+1,rd) \to \mathbb{P}(n+1,rd) $ is a $GL_{n+1}/\mu_{rd}$ equivariant morphism. Consider the normal subgroup $\mathbb{G}_m/\mu_{rd}$ embedded diagonally in $GL_{n+1}/\mu_{rd}$. Note that $(GL_{n+1}/\mu_{rd})/(\mathbb{G}_m/\mu_{rd}) \cong PGL_{n+1}$ and that $f: \mathbb{A}(n+1,rd) \to \mathbb{P}(n+1,rd) $ is a  $\mathbb{G}_m/\mu_{rd}$ torsor. Then the conclusion follows from \cite{KS}, Lemma $2.3$.\QEDB



\begin{corollary}\label{tautological family main}
Both $\mathscr{H}_{n,r,d}$ and $M_{n,r,d}$ are unirational. Moreover,
\begin{itemize}
\item[(1)] if $\textrm{gcd}(rd, n+1) \nmid d$, the stack $\mathscr{H}_{n,r,d}$ is not rational,
\item[(2)] if $\textrm{gcd}(rd, n+1) \mid d$, then $\mathscr{H}_{n,r,d}$ is rational if and only if $M_{n,r,d}$ is rational.

\end{itemize}

\end{corollary}

\noindent\textbf{Proof.} Unirationality of both $\mathscr{H}_{n,r,d}$ and $M_{n,r,d}$ are clear from the quotient structure as seen in Lemma ~\ref{Coarse moduli of plane curves 1}.

Note that if a \textcolor{black}{Deligne-Mumford} stack $\mathscr{X}$ \textcolor{black}{with finite inertia} is rational then its coarse moduli space $M$ is rational and the stack admits a section from a Zariski open subset of the coarse moduli. Conversely, if $\mathscr{X}$ is a gerbe over a Zariski open subset of $M$, then the rationality of $M$ and the existence of a section from a Zariski open subset of $M$ implies the rationality of the stack $\mathscr{X}$. Hence $\mathscr{H}_{n,r,d}$ is rational if and only if $M_{n,r,d}$ is rational and $\mathscr{H}_{n,r,d}$ admits a section from a  Zariski open subset of $M_{n,r,d}$. The conclusion follows from Theorem \ref{tautological family}.\QEDB

\begin{corollary}\label{rath1}
Assume $rd\geq 4$.
 \begin{itemize}
        \item[(i)] If $rd$ is odd, then $\mathscr{H}_{1,r,d}$ is rational. 
        \item[(ii)] If $rd$ is even then $\mathscr{H}_{1,r,d}$ is rational if and only if $d$ is even.
    \end{itemize}
\end{corollary}
\noindent\textbf{Proof.}  The proof follows from Corollary ~\ref{tautological family main} and the fact that the moduli of points on $\mathbb{P}^1$ denoted by $M_{0,n}$ is rational for $n \geq 4$.\QEDB

\begin{corollary}\label{rath2}
Assume $rd\geq 4$.
\begin{itemize}
     \item[(i)] If $rd \geq 49$ the stack $\mathscr{H}_{2,r,d}$ is rational if and only if either $3 \mid d$ or $3 \nmid rd$.
     \item[(ii)] If $(r,rd) \in S:=\left\{(3,6), (3,12), (3,15), (9,18), (3,24), (6,24), (12,24), (3,48), (6,48)\right\}$, then $\mathscr{H}_{2,r,d}$ is not rational.
     \item[(iii)] If $rd \leq 48$ and $(r,rd) \notin S$, $\mathscr{H}_{2,r,d}$ is rational if and only if $C(2,rd)$ is rational.
\end{itemize}
\end{corollary}
\noindent\textbf{Proof.} Follows from Theorem \ref{tautological family}, Corollary ~\ref{tautological family main} and results on rationality of plane curves in \cite{BB10} and \cite{BBK09}. \hfill\QEDB

\subsubsection{Picard group \textcolor{black}{of $M_{n,r,d}^0$} and \textcolor{black}{non} existence of a tautological family over $M_{n,r,d}^0$} 

\textcolor{black}{In this subsection we compute the Picard group of $M_{n,r,d}^0$ for $n = 1,2$ by showing that in those cases the codimension of the complement of $M_{n,r,d}^0$ in $M_{n,r,d}$ is at least two, which is a condition required for Theorem ~\ref{tautological family}, $(3)$. The case $n = 2$ is the major case where we show that plane curves with non-trivial linear automorphisms have codimension at least two in $M_{2,r,d}$ using results of \cite{BB}. In proposition ~\ref{tautological family in any dimension} we} show that there does not exist a tautological family on $M_{n,r,d}^0$ \textcolor{black}{for any $n$}.

\begin{lemma}\label{comparison of picard groups}
Suppose that $k \mid r$. \textcolor{black}{Suppose that $GL_{n+1}/\mu_{kd}$ acts on $\mathbb{A}_{sm}(n, rd)$ with the action given by $[A]\cdot f(x)=f(A^{-1}x)$.  }
\begin{itemize}
    \item[(1)] The natural map $$\textrm{Pic}([\mathbb{A}_{sm}(n, rd) / GL_{n+1}/\mu_{kd}]) \to \textrm{Pic}(\mathscr{H}_{n,r,d})$$ is injective of index $$\frac{k \cdot \textrm{gcd}(d,n+1)}{\textrm{gcd}(kd,n+1)}.$$
    \item[(2)] $\textrm{Pic}(\mathscr{H}_{n,r,d})$ is generated by $$\textrm{det}^{\frac{d}{\textrm{gcd}(d,n+1)}}$$ 
    \textrm{satisfying} $\textrm{det}^{rd(rd-1)^n} = \textrm{id}$ (see \cite{AV}, Theorem $5.1$)
    \item[(3)] $\textrm{Pic}([\mathbb{A}_{sm}(n, rd) / GL_{n+1}/\mu_{kd}])$ is generated by $$\textrm{det}^{\frac{kd}{\textrm{gcd}(kd,n+1)}}$$ satisfying $\textrm{det}^{rd(rd-1)^n} = \textrm{id}$. \textcolor{black}{Setting $r = k$, one obtains $\textrm{Pic}(\mathscr{D}_{n,r,d})$ is generated by $$\textrm{det}^{\frac{rd}{\textrm{gcd}(rd,n+1)}}$$ satisfying $\textrm{det}^{rd(rd-1)^n} = \textrm{id}$.}

\end{itemize}
\end{lemma}

\noindent\textbf{Proof.} We first note the following chain of maps of algebraic groups
$$\mathbb{G}_m \hookrightarrow GL_{n+1} \rightarrow GL_{n+1}/\mu_d \rightarrow GL_{n+1}/\mu_{kd}.$$
This induces the following chain of inclusion of character groups
$$\widehat{GL_{n+1}/\mu_{kd}} \hookrightarrow \widehat{GL_{n+1}/\mu_{d}} \hookrightarrow \widehat{GL_{n+1}} \hookrightarrow \widehat{\mathbb{G}_m}.$$
The fact that the first two maps of character groups is injective follows from the fact that $$GL_{n+1} \to GL_{n+1}/\mu_{d} \to GL_{n+1}/\mu_{kd}$$ are surjective maps. The fact that the last map is injective is due to the fact that the map is induced by the restriction of the character $\textrm{det}$ of ${GL_{n+1}}$ to $\mathbb{G}_m \hookrightarrow GL_{n+1}$ embedded diagonally. Let us note the index of the injective group homomorphisms. The index of $\widehat{GL_{n+1}} \hookrightarrow \widehat{\mathbb{G}_m }$ is $n+1$. The index of $\widehat{GL_{n+1}/\mu_d} \hookrightarrow \widehat{GL_{n+1}}$ is $\frac{d}{\textrm{gcd}(d,n+1)}$ (see \cite{AV}, Proof of Theorem $5.1$). Similarly the index of $\widehat{GL_{n+1}/\mu_{kd}} \hookrightarrow \widehat{GL_{n+1}}$ is $\frac{kd}{\textrm{gcd}(kd,n+1)}$. Hence index of $\widehat{GL_{n+1} /\mu_{kd}} \hookrightarrow \widehat{GL_{n+1} /\mu_{d}}$ is $\frac{k \cdot \textrm{gcd}(d,n+1)}{\textrm{gcd}(kd,n+1)}$. 

Now let $V = \mathbb{A}(n,rd)$, $X = \mathbb{A}_{sm}(n,rd)$ so that $\Delta = V\setminus X$ is the irreducible locus of singular forms. Using \cite{EG}, Proposition $18$ and the following exact sequences ( see \cite{AV}, Proof of Theorem $5.1$) 
$$A_{n-1}^{GL_{n+1}/\mu_d}(V\setminus X) \xrightarrow{f} \widehat{GL_{n+1}/\mu_d} \to \textrm{Pic}(\mathscr{H}_{n,r,d}) \to 0,$$
$$A_{n-1}^{GL_{n+1}/\mu_{kd}}(V\setminus X) \xrightarrow{g} \widehat{GL_{n+1}/\mu_{kd}} \to \textrm{Pic}([\mathbb{A}_{sm}(n,rd)/\mu_{kd}]) \to 0,$$
we have that the image of the map $f$ or $g$ is the subgroup generated by $\Delta$ \textcolor{black}{and we have the} following sequence of injective group homomorphisms
$$\Delta \hookrightarrow \widehat{GL_{n+1}/\mu_{kd}} \hookrightarrow \widehat{GL_{n+1}/\mu_{d}} \hookrightarrow \widehat{GL_{n+1}} \hookrightarrow \widehat{\mathbb{G}_m}.$$
Hence we have the following commutative diagram where the rows are exact.
\[
\begin{tikzcd}
0 \arrow[r] & \Delta \arrow[r] \arrow[equal]{d}  & \widehat{GL_{n+1}/\mu_{d}} \arrow[r] & \textrm{Pic}(\mathscr{H}_{n,r,d}) \arrow[r] & 0 \\
0 \arrow[r] & \Delta \arrow[r]  & \widehat{GL_{n+1}/\mu_{kd}} \arrow[r] \arrow[u] & \textrm{Pic}([\mathbb{A}_{sm}(n,rd)/\mu_{kd}]) \arrow[r] \arrow[u] & 0
\end{tikzcd}
\]
This implies that the index of $\textrm{Pic}([\mathbb{A}_{sm}(n,rd)/GL_{n+1}/\mu_{kd}])$ in $\textrm{Pic}(\mathscr{H}_{n,r,d})$ is equal to the index of $\widehat{GL_{n+1}/\mu_{kd}}$ in $\widehat{GL_{n+1}/\mu_{d}}$ which is equal to 
$\frac{k \cdot \textrm{gcd}(d,n+1)}{\textrm{gcd}(kd,n+1)}$. This proves $(1)$.

 The index of $\Delta$ in $\mathbb{G}_m$ for the given action is $rd(n+1)(rd-1)^n$. This together with the index of $\widehat{GL_{n+1}/\mu_{kd}}$ in $\widehat{\mathbb{G}_m}$ gives $(2)$ and $(3)$.\QEDB 

\vspace{5pt}

\textcolor{black}{In the following two propositions we prove that the codimension of the complement of $M_{n,r,d}^0$ in $M_{n,r,d}$ is at least two for $n = 1,2$.}
 

\begin{proposition}\label{codimension of stacky locus of orbifold for $n = 1$}
The complement of $M_{1,r,d}^0$ in $M_{1,r,d}$ has codimension $\geq 2$ for $rd \geq 8$. \textcolor{black}{Hence $\textrm{Pic}(M_{1,r,d}^0) = \textrm{Pic}(\mathscr{D}_{1,r,d}^0) = \textrm{Pic}(\mathscr{D}_{1,r,d})$.}
\end{proposition}

\noindent\textbf{Proof.} The proof of the above theorem is a verbatim consequence of the methods used in \cite{GV2}, Proposition $2.1$ and hence we omit the proof. \textcolor{black}{The assertion on Picard groups then follow from the fact that $\mathscr{M}_{1,r,d}^0 \cong \mathscr{D}_{1,r,d}^0$ and $\mathscr{D}_{1,r,d}$ is smooth.}\QEDB

\vspace{5pt}

We aim to prove an analogue of the \textcolor{black}{Proposition ~\ref{codimension of stacky locus of orbifold for $n = 1$}} above in dimension two. Using a stratification of the locus plane curves of degree $k$ with non-trivial linear automorphisms, we can show that analogous conclusion holds for $M_{2,r,d}$.

\begin{proposition}\label{codimension of stacky locus of orbifold for $n = 2$}
If $\textrm{char}(\mathbb{k}) = 0$ or $\textrm{char}(\mathbb{k}) > (rd-1)(rd-2) + 1$, the complement of $M_{2,r,d}^0$ in $M_{2,r,d}$ has codimension greater than or equal to $2$ for $rd \geq 7$. \textcolor{black}{Hence $\textrm{Pic}(M_{2,r,d}^0) = \textrm{Pic}(\mathscr{D}_{2,r,d}^0) = \textrm{Pic}(\mathscr{D}_{2,r,d})$.}
\end{proposition}

\noindent\textbf{Proof.} We use the stratification given in \cite{BB}, Theorem $7$ of the locus of plane curves with non-trivial cyclic automorphisms into six stratas, to give upper bounds to $\rho_{m,a,b}(M_g^{Pl}(\mathbb{Z}/m\mathbb{Z}))$ which is the loci of smooth plane curves of degree $rd$ containing the subgroup of $PGL_3$, abstractly isomorphic to $\mathbb{Z}/m\mathbb{Z}$ and which injects into $PGL_3$ by the map $\zeta \to \textrm{diag}\langle 1, \zeta^a,  \zeta^b \rangle$ and $m$ is a prime. Let $L_{j,*}$ denote a homogeneous polynomial of degree $j$ such that the variable $* \in \{X,Y,Z\}$ does not appear. In what follows, we follow the notations of \cite{BB}, Theorem 7.

(1) $ C \in \rho_{m,0,1}(M_g^{Pl}(\mathbb{Z}/m\mathbb{Z}))$ with $m \mid rd-1$ and any such polynomial is of the form

\begin{equation*}
   Z^{rd-1}Y +  \sum_{j \in S(2)_m} Z^{rd-j}L_{j,Z} + L_{rd,Z}
\end{equation*}
where $S(2)_m = \left\{2 \leq j \leq rd-1 | m \mid rd-j \right\}$.
Let $rd-1 = mk$. Then 
$$\textrm{dim}\left(\rho_{m,a,b}(M_g^{Pl}(\mathbb{Z}/m\mathbb{Z}))\right) = \left(\sum_{j \in S(2)_m} (j+1)\right) + (rd+1)$$
$$\implies\textrm{dim}\left(\rho_{m,a,b}(M_g^{Pl}(\mathbb{Z}/m\mathbb{Z}))\right)  = \left(\sum_{j \in S(2)_m} j\right) + k + rd  = \frac{mk(k-1)}{2} + 2k + rd. $$ 

The dimension of $\textrm{C}(rd)$ the moduli space of plane curves of degree $rd$ is $\frac{(rd+2)(rd+1)}{2} - 8$. Hence we want to see when 
\begin{equation}\label{eq1}
    \frac{mk(k-1)}{2} + 2k + rd > \frac{(rd+2)(rd+1)}{2} - 10.
\end{equation}
The inequality \eqref{eq1} implies 
\begin{equation}\label{eq2}
    (m-1)(rd)^2 + rd(2m-2) - (19m-3) < 0.
\end{equation}
Using $rd \geq m + 1$, we have 
$m(m^2 + 3m -20) < 0$ which gives $m = 2$ or $m = 3$.
Plugging the value of $m$ in \eqref{eq2}, we have $(rd)^2 + 2(rd) < 35$. This gives $rd = 5$ considering that $m \mid rd-1$. Plugging in $m = 3$ we have that $rd = 4$. Hence for $rd \geq 6$, $\textrm{dim}\left(\rho_{m,a,b}(M_g^{Pl}(\mathbb{Z}/m\mathbb{Z}))\right)$ has codimension greater than or equal to $2$ in $\textrm{C}(rd)$.

(2) $ C \in \rho_{m,0,1}(M_g^{Pl}(\mathbb{Z}/m\mathbb{Z}))$ with $m \mid rd$ and any such polynomial is of the form

\begin{equation*}
   Z^{rd} +  \sum_{j \in S(1)_m} Z^{rd-j}L_{j,Z} + L_{rd,Z}
\end{equation*}
where $S(1)_m = \left\{1 \leq j \leq rd-1 | m \mid rd-j \right\}$.
Using identical computation we have that 
$\textrm{dim}(\rho_{m,a,b}(M_g^{Pl}(\mathbb{Z}/m\mathbb{Z})))$ has codimension greater than or equal to $2$ in $\textrm{C}(rd)$ unless $(rd,m) = (4,2)$ or $(6,2)$. Hence in this case $$\textrm{dim}\left(\rho_{m,a,b}(M_g^{Pl}(\mathbb{Z}/m\mathbb{Z}))\right)$$ has codimension greater than or equal to $2$ if $rd \geq 7$.

(3) $ C \in \rho_{m,a,b}(M_g^{Pl}(\mathbb{Z}/m\mathbb{Z}))$ with $m \mid (rd)^2-3rd+2$, $(a,b) \in \Gamma_m$ where
\begin{equation*}
    \Gamma_m = \{(a,b) \in \mathbb{N}^2 \mid \textrm{gcd}(a,b) = 1, 1 \leq a \neq b \leq m-1 \}
\end{equation*}
and any such polynomial is of the form
\begin{equation*}
 X^{rd-1}Y + Y^{rd-1}Z + Z^{rd-1}X + \sum_{j=2}^{rd-2} \left (\sum_{j \in S(1)^{j,X}_{m,(a,b)}} \beta_{j,i} X^{rd-i} Y^i Z^{j-i} + \sum_{j \in S^{j,Y}_{m,(a,b)}} \alpha_{j,i} X^{j-i} Y^{rd-j} Z^{i} + \sum_{j \in S^{j,Z}_{m,(a,b)}} \gamma_{j,i} X^{j-i} Y^i Z^{rd-j} \right) 
\end{equation*}
where $S(1)^{j,X}_{m,(a,b)} = \left\{1 \leq i \leq j \mid ai + b(j-i) \equiv a (\textrm{mod}(m)) \right\}$, $S^{j,Y}_{m,(a,b)} = \left\{1 \leq i \leq j \mid bi + a(rd-j) \equiv a (\textrm{mod}(m)) \right\}$ and $S^{j,Z}_{m,(a,b)} = \left\{1 \leq i \leq j \mid ai + b(rd-j) \equiv a (\textrm{mod}(m)) \right\}$
Note that we have, $a$, $b$ and $a-b$ are coprime to $m$ since all of these are strictly less than $m$. Hence 
$$\textrm{dim}\left(\rho_{m,a,b}(M_g^{Pl}(\mathbb{Z}/m\mathbb{Z}))\right) \leq \sum\limits_{j=2}^{rd-2} \frac{3}{m}(j+1). $$ 
So we want to see when
\begin{equation}\label{eq3}
    \frac{3}{m}\left[\frac{(rd-1)(rd-2)}{2} + (rd-4)\right] > \frac{(rd+2)(rd+1)}{2} - 10.
\end{equation}
Since $m \geq 3$ (otherwise $\Gamma_m$ is empty), \eqref{eq3} implies that 
\begin{equation}\label{eq4}
    \frac{(rd-2)(rd-1)}{2} + (rd-4) > \frac{(rd+2)(rd+1)}{2} - 10.
\end{equation}
This implies $rd < 3$.

(4) In all of these cases, noting that $a-b$ and $m$ are coprime, and that $m \geq 3$, one can see that 
$$\textrm{dim}\left(\rho_{m,a,b}(M_g^{Pl}(\mathbb{Z}/m\mathbb{Z}))\right) \leq \frac{1}{3}\left(\left(\sum_{j = 2}^{rd-2}(j+1)\right) + 3rd-8\right).$$
Hence we want to see when
$$2\left(\sum_{j = 2}^{rd-2}(j+1) + 3rd-8 \right) > 3(rd)^2 + 9rd -54.$$
But that implies $2(rd)^2 + 6rd -40 < 0$ which means $rd \leq 3$.

$(5), (6)$ We need $rd \leq 3$ and this follows from identical computation as in $(4)$. \newline
\textcolor{black}{The assertion on Picard groups then follow from the fact that $\mathscr{M}_{2,r,d}^0 \cong \mathscr{D}_{2,r,d}^0$ and $\mathscr{D}_{2,r,d}$ is smooth.}\QEDB

\vspace{5pt}

\textcolor{black}{In the following proposition we show that even if there is a tautological family over some open subscheme of $M_{n,r,d}$, that open subscheme can never be taken to be $M_{n,r,d}^0$}

\begin{proposition}\label{tautological family in any dimension}
There does not exist a \textcolor{black}{a tautological family on $M_{n,r,d}^0$ or in other words there does not exist a} global section $s: M_{n,r,d}^0 \to \mathscr{H}_{n,r,d}$ of the coarse moduli map $p: \mathscr{H}_{n,r,d} \to M_{n,r,d}$. 
\end{proposition}

\color{black}
\noindent\textbf{Proof.} We use the notations of Lemma ~\ref{relative picard of universal Brauer Severi} and its proof. Recall that we have shown in Lemma ~\ref{relative picard of universal Brauer Severi} that $\textrm{Pic}(P) = \textrm{Pic}^{GL_{n+1}/\mu_{rd}}(\mathbb{A}_{sm}^0(n,rd) \times \mathbb{P}^n)$ where $(P \to M_{n,r,d}^0, D)$ is the pair consisting of the universal Brauer-Severi scheme and the universal divisor of degree $rd$ on each fibre given by the isomorphism $M_{n,r,d}^0 \cong \mathscr{D}_{n,r,d}^0$. If there is a tautological family on $M_{n,r,d}^0$, then $D$ must be $r$-divisible in $\textrm{Pic}(P)$, or in other words the $GL_{n+1}/\mu_{rd}$ invariant incidence divisor $\mathbb{D}$ on $\mathbb{A}_{sm}^0(n,rd) \times \mathbb{P}^n$ must be $r$- divisible in $\textrm{Pic}^{GL_{n+1}/\mu_{rd}}$$(\mathbb{A}_{sm}^0(n,rd) \times \mathbb{P}^n)$. To show that this cannot happen we first claim that $$\textrm{Pic}^{GL_{n+1}/\mu_{rd}}(\mathbb{A}_{sm}^0(n,rd) \times \mathbb{P}^n) = \textrm{Pic}^{PGL_{n+1}}(\mathbb{B}_{sm}^0(n,rd) \times \mathbb{P}^n)$$ where $\mathbb{B}_{sm}^0(n,rd)$ is the projectivization of $\mathbb{A}_{sm}^0(n,rd)$. We use Lemma $2.3$ of \cite{KS}. Note that there is a (free) diagonal action of $GL_{n+1}/\mu_{rd}$ on both $\mathbb{A}_{sm}^0(n,rd) \times \mathbb{P}^n$ and $\mathbb{B}_{sm}^0(n,rd) \times \mathbb{P}^n$. Consider the natural map from $$\mathbb{A}_{sm}^0(n,rd) \times \mathbb{P}^n \to \mathbb{B}_{sm}^0(n,rd) \times \mathbb{P}^n$$ which is identity on $\mathbb{P}^n$. This map is equivariant under the action of $GL_{n+1}/\mu_{rd}$. Note that the induced action of the closed normal subgroup $\mathbb{G}_m/\mu_{rd} \hookrightarrow GL_{n+1}/\mu_{rd}$ as classes of scalar matrices preserves the fibres of the equivariant map as $$ [\beta]\cdot (\lambda f, p) = (\lambda f(\beta^{-1}x), [\beta]p = p) $$ We show that this action on the fibres is free and transitive. We first show the freeness of the action. Let $([f],p) \in \mathbb{B}_{sm}^0(n,rd) \times \mathbb{P}^n$. Then a point on the fibre is of the form $(\lambda f, p)$. Then $[\beta] \in \mathbb{G}_m/\mu_{rd}$ acts on $(\lambda f, p)$ to give $(\lambda f(\beta^{-1}x), p) = (\beta^{-rd} \cdot \lambda f(x), p)$. If $[\beta]$ is in the stabilizer of $(\lambda f, p)$, then $[\beta]$ is in the stabilizer of $\lambda f$ (since the action is diagonal) and then $\beta^{-rd} = 1$ and hence $[\beta] = \textrm{Id}_{\mathbb{G}_m/\mu_{rd}}$. This establishes that $\mathbb{G}_m/\mu_{rd}$ acts freely on the fibres. To show transitivity, note that for any $\lambda \in \mathbb{G}_m$, if we choose $\beta$ to be an $r$-th root of $\lambda$, then $[\beta] \cdot (f,p) = (\lambda f, p)$. Hence the  
the natural $GL_{n+1}/\mu_{rd}$ equivariant map from $$\mathbb{A}_{sm}^0(n,rd) \times \mathbb{P}^n \to \mathbb{B}_{sm}^0(n,rd) \times \mathbb{P}^n$$ is a $\mathbb{G}_m/\mu_{rd}$ torsor.
Now considering that $\displaystyle\frac{GL_{n+1}/\mu_{rd}}{\mathbb{G}_m/\mu_{rd}} \cong PGL_{n+1}$, we have by Lemma $2.3$, \cite{KS} an isomorphism of quotient stacks 
$$[\displaystyle\frac{\mathbb{A}_{sm}^0(n,rd) \times \mathbb{P}^n}{GL_{n+1}/\mu_{rd}}] \cong [\displaystyle\frac{\mathbb{B}_{sm}^0(n,rd) \times \mathbb{P}^n}{PGL_{n+1}}]$$
Therefore comparing their Picard groups we have that $$\textrm{Pic}^{GL_{n+1}/\mu_{rd}}(\mathbb{A}_{sm}^0(n,rd) \times \mathbb{P}^n) = \textrm{Pic}^{PGL_{n+1}}(\mathbb{B}_{sm}^0(n,rd) \times \mathbb{P}^n)$$
Now consider the image of the incidence divisor $\mathbb{D}$ (which we still call $\mathbb{D}$ by a slight abuse of notation) inside $\textrm{Pic}(\mathbb{B}_{sm}^0(n,rd) \times \mathbb{P}^n) = \textrm{Pic}(\mathbb{B}_{sm}^0(n,rd) \times \textrm{Pic}(\mathbb{P}^n)$. We see that its class is given by $(\mathscr{O}_{\mathbb{B}_{sm}^0(n,rd)}(1), \mathscr{O}_{\mathbb{P}^n}(rd))$. This is not $r$-divisible in
$\textrm{Pic}(\mathbb{B}_{sm}^0(n,rd) \times \mathbb{P}^n)$ and is hence not $r$-divisible in $\textrm{Pic}^{PGL_{n+1}}(\mathbb{B}_{sm}^0(n,rd) \times \mathbb{P}^n)$. \QEDB
\color{black}

\vspace{5pt}

\subsubsection{Picard group of coarse moduli of \textcolor{black}{the stack of} cyclic covers over the projective space}

Now we prove the final result of this section where we compute the Picard group of the coarse moduli.
\begin{theorem}\label{Picard Group of the coarse moduli}
  Pic$(M_{n, r, d}) = 0$ if $rd \geq 4$. 
\end{theorem}

\noindent\textbf{Proof.} Let $g = \textrm{gcd}(rd, n+1)$. Note that by Lemma ~\ref{comparison of picard groups}, $\textrm{Pic}(\mathscr{D}_{n, r, d})$ is generated by $s = \textrm{det}^{rd/g}$ satisfying the relation $s^{(rd-1)^n g} = 1$. Now by \cite{KKV}, Proposition $4.2$, $\textrm{Pic}( M_{n,r,d})$ injects inside $\textrm{Pic}(\mathscr{D}_{n,r,d})$  and is contained inside the subgroup generated by elements which restricted to any non-trivial stabilizer group are the image of the identity character. Since this is also a finite cyclic group we have that it is generated by an element $s_{\alpha} = s^{\alpha}$. Now since $rd \geq 4$ by \cite{ZZ} Lemma $2.3$ we have that there exist a smooth curve with an element of order $(rd-1)^n$ in its stabilizer group with presentation $\rho =  \textrm{diag} \langle \zeta, \zeta^{(1-rd)}, \zeta^{(1-rd)^2}, ... , \zeta^{(1-rd)^n} \rangle$ where $\zeta$ is \textbf{} $(rd-1)^n$-th root of unity. Note that $m = 1 + \sum_{i = 1}^{n} (1-rd)^i$ is coprime to $(rd-1)^n$. Now we have that $s_{\alpha}(\rho) = 1$ and hence $\zeta^{mrd \alpha/g} = 1 \implies (rd - 1)^n \mid mrd \alpha/g$. Since $mrd$ is coprime to $(rd - 1)^n$, we have that $(rd - 1)^n \mid \alpha$. \par 
We now show that $g$ divides $\alpha$. Note that the polynomial $\sum_{i = 0}^{n} x_i^{rd}$ is a smooth polynomial with an element order $rd$ in its stabilizer group with representation diag $\langle 1,..., \zeta_{rd} \rangle$ where $\zeta_{rd}$ is a \textbf{} $rd$- th root of unity. Then $\zeta_{rd}^{rd \alpha/g} = 1 \implies rd \mid rd \alpha/ g \implies g \mid \alpha$. \par
Since $g \mid rd$ we have that $\textrm{gcd}(g, (rd-1)^n) = 1$. Hence $(rd-1)^ng \mid \alpha.$ So $s_{\alpha} = 1$. Hence $\textrm{Pic}(M_{n, r,d}) = 0$.\QEDB

\section{Stacks of cyclic triple covers over \texorpdfstring{$\mathbb{P}^{1}$}{}}\label{3}

Throughout this section, we work with the stack of cyclic triple covers of $\mathbb{P}^1$ with the objective of proving Theorem ~\ref{main2intro}.

\subsection{Rigidification of stacks of cyclic triple covers over \texorpdfstring{$\mathbb{P}^1$}{}} We formally introduce the main object of our study, namely, the fibred category $\mathscr{H}_{1,3,d_1,d_2}$. We need to fix some notations first that we will follow throughout this section. We start with positive integers $d_{1},d_{2}$ such that $0\leq d_{1}\leq 2d_{2},0\leq d_{2}\leq 2d_{1}$. Set
$$l_{1}=2d_{1}-d_{2},\quad l_{2}=2d_{2}-d_{1}.$$
If $l_i=0$ for any $i$ then the stack is isomorphic to the stack of simple cyclic covers that has been treated in the previous section. We assume throughout the rest of this section that $l_i\geq 4$ for $i=1,2$ and assume that $\textrm{char}(\mathbb{k})$ does not divide $2l_{1}l_{2}$.
We define the group $\Gamma(d_1,d_2)$ as follows \textcolor{black}{(see \cite{AV}, Theorem $6.5$)} $$\Gamma(d_{1},d_{2}):=G_{m}\times GL_{2}/(\mu_{d_{1}}\times \mu_{d_{2}}) \cong \textbf{\underline{Aut}}(\mathbb{P}_{\mathbb{Z}}^n, \mathscr{O}(-d_{1}))\times_{PGL_{2}} \textbf{\underline{Aut}}(\mathbb{P}_{\mathbb{Z}}^n, \mathscr{O}(-d_{2})).$$ 
Here the we denote the quotient with respect to the embedding
$$\mu_{d_{1}}\times \mu_{d_{2}}\xhookrightarrow{i} G_{m}\times GL_{2}$$
$$(x_{1},x_{2})\rightarrow (x_{2}/x_{1},\ x_{1} I_{2\times 2})$$ where $I_{2\times 2}$ denotes the identity in $GL_{2}$ and $\mu_{d_{i}}$ denote the $d_{i}$-th root of unity in $G_{m}.$

\smallskip

We now introduce the main object of our study, namely the stack $\mathscr{H}_{1,3,d_{1},d_{2}}$.
Recall that if $\pi: X \to Y$ be an abelian cover with group $\mu_3$ (i.e., $Y \cong X/\mu_3$) with $X$ normal and $Y$ smooth then by \cite{P91}, Theorem $2.1$, $X$ is uniquely determined by line bundles $L_1$ and $L_2$ and divisors $D_1 \in |L_1^{\otimes -2} \otimes L_2|$ and $D_2 \in |L_2^{\otimes -2} \otimes L_1|$. If further $X$ and $Y$ are smooth projective curves, then by \cite{P91}, Theorem $3.1$, both $D_1$ and $D_2$ are smooth and disjoint.  

\begin{definition}\label{Definition of stack H}
Let $\mathscr{H}_{1,3,d_{1},d_{2}}$ denote the fibred category 
\color{black}
\begin{itemize}
    \item[(a)]  whose objects over any $\mathbb{k}$-scheme $S$  consist of a Brauer--Severi scheme $P \rightarrow S$ of relative dimension 1, line bundles $\mathscr{L}_{1},\mathscr{L}_{2}$ on $P$ which restrict to degree $-d_{1},-d_{2}$ on the fibre of every geometric point $s$ of $S$ respectively and injective homomorphisms $i_{1} : \mathscr{L}_{1}^{\otimes 2} \rightarrow \mathscr{L}_{2}$ and $i_{2} : \mathscr{L}_{2}^{\otimes 2} \rightarrow \mathscr{L}_{1}$ that remains injective with smooth disjoint associated divisors when restricted to any geometric fibre
    \item[(b)] and morphisms consist of Cartesian diagrams between Brauer--Severi schemes preserving the line bundles and the injective homomorphisms. 
\end{itemize}
\color{black}
\end{definition}

It has been shown in \cite{AV}, Section 6 that this stack is isomorphic to the stack of smooth cyclic triple covers over $\mathbb{P}^1$.
Consider the following open subscheme $$\mathbb{U}:=\mathbb{A}_{sm}(1,l_1)\times \mathbb{A}_{sm}(1,l_2)- Z$$ where $Z$ is the closed subscheme consisting of pairs of forms with common root. Equivalently, $\mathbb{U}$ is the open subscheme consisting of pairs of smooth forms with no common root.

\begin{theorem}\label{avtriple}(\cite{AV}, Theorem 6.5)
$\mathscr{H}_{1,3,d_{1},d_{2}} \cong [\mathbb{U}/\Gamma(d_{1},d_{2})]$
under the action $\Gamma(d_{1},d_{2}) \times \mathbb{U}\rightarrow \mathbb{U}$ given by
$$[(\alpha,A)]\cdot (f_{1}(x),f_{2}(x)) =  (\alpha)^{d_{2}} \cdot (f_{1}(A^{-1}x),(\alpha)^{-2d_{2}} \cdot f_{2}(A^{-1}x))$$
where $(\alpha,A)\in G_{m}\times GL_{2} \And [(\alpha,A)]\in \Gamma(d_{1},d_{2})$ denotes the corresponding class under the quotient by $i$.
\end{theorem}

Next we introduce the stack $\mathscr{D}_{1,3,d_{1},d_{2}}$ whose objects and morphisms are introduced below. 

\begin{definition}\label{definition of orbifold D}
Let $\mathscr{D}_{1,3,d_{1},d_{2}}$ denote the fibred category 
\color{black}
\begin{itemize}
    \item[(a)]  whose objects over any $\mathbb{k}$-scheme $S$ consists of a Brauer--Severi scheme $P \rightarrow S$ of relative dimension $1$, line bundles $\mathscr{L}_{1},\mathscr{L}_{2}$ on $P$ which restrict to degree $-l_{1},-l_{2}$ \textcolor{black}{(recall that $l_1 = 2d_1-d_2$ and $l_2 = 2d_2-d_1$)} on a fibre of every geometric point $s$ of $S$ respectively and injective homomorphisms $i_{1} : \mathscr{L}_{1} \rightarrow \mathscr{O}_{P}$ and $i_{2} : \mathscr{L}_{2} \rightarrow \mathscr{O}_{P}$, that remains injective with smooth disjoint associated divisors when restricted to any geometric fibre
    \item[(b)] and morphisms consist of Cartesian diagrams between Brauer--Severi schemes preserving the line bundles and the injective homomorphisms. 
\end{itemize}
\color{black}

\end{definition}

We give a quotient stack structure to the stack $\mathscr{D}_{1,3,d_{1},d_{2}}$ following the treatment of \cite{AV}.

\begin{theorem}\label{Quotient Stack Structure}
$\mathscr{D}_{1,3,d_{1},d_{2}} \cong [\mathbb{U}/\Gamma(l_{1},l_{2})]$ under the action $\Gamma(l_{1},l_{2})\times \mathbb{U}\rightarrow \mathbb{U}$ given by
$$[(\alpha,A)]\cdot (f_{1}(x),f_{2}(x))\rightarrow (f_{1}(A^{-1}x),(\alpha)^{-l_{2}}\cdot f_{2}(A^{-1}x)),$$
where $(\alpha,A)\in G_{m}\times GL_{2}$ and $[(\alpha,A)]\in \Gamma(l_{1},l_{2})$ denotes the corresponding class under the quotient by $i$.
\end{theorem}
\noindent\textbf{Proof.} Consider the auxiliary fibred category $\widetilde{\mathscr{D}}_{1,3,d_{1},d_{2}}$ whose objects over a $\mathbb{k}$-scheme $S$ consist of an object $(P \to S, \mathscr{L}_{1},\mathscr{L}_{2}, i_{1}: \mathscr{L}_{1} \to \mathscr{O}_P, i_{2}: \mathscr{L}_{2} \to \mathscr{O}_P)$ in $\mathscr{D}_{1,3,d_{1},d_{2}}$ along with an isomorphism $\phi_{1}: (P, \mathscr{L}_{1}) \cong (\mathbb{P}_S^1, \mathscr{O}_{\mathbb{P}_S^1}(-l_{1}))$ and $\phi_{2}: (P,\mathscr{L}_{2}) \cong (\mathbb{P}_S^1, \mathscr{O}_{\mathbb{P}_S^1}(-l_{2}))$ such that $\phi_{1}\circ \phi_{2}^{-1}$ induces identity from $\mathbb{P}^1_{S} \to \mathbb{P}^1_{S}$. Note that there is forgetful map from $\widetilde{\mathscr{D}}_{1,3,d_{1},d_{2}} \to \mathscr{D}_{1,3,d_{1},d_{2}}$. We define a base preserving functor from $\widetilde{\mathscr{D}}_{1,3,d_{1},d_{2}} \to \mathbb{U}$. For any object of $\widetilde{\mathscr{D}}_{1,3,d_{1},d_{2}}(S)$ consider the composite homomorphism $\phi_{j} \circ i_{j} \circ \phi^{-1}_{j}: \mathscr{O}_{\mathbb{P}_S^1}(-l_{1}) \to \mathscr{O}_{\mathbb{P}_S^1}, j=1,2$ which gives an element of $\mathbb{U}(S)$. Conversely given an element of $\mathbb{U}(S)$, it can be viewed as an injective hommorphism $j_{1}: \mathscr{O}_{\mathbb{P}_S^1}(-l_{1}) \to \mathscr{O}_{\mathbb{P}_S^1},j_{2}:\mathscr{O}_{\mathbb{P}_S^1}(-l_{2}) \to \mathscr{O}_{\mathbb{P}_S^1}$ which remains injective with smooth associated divisor when restricted to any geometric fibre. Then consider the element $(\mathbb{P}_S^1 \to S, \mathscr{O}_{\mathbb{P}_S^1}(-l_{1}), j_{1}: \mathscr{O}_{\mathbb{P}_S^1}(-l_{1}) \to \mathscr{O}_{\mathbb{P}_S^1},\mathscr{O}_{\mathbb{P}_S^1}(-l_{2}), j_{2}: \mathscr{O}_{\mathbb{P}_S^1}(-l_{2}) \to \mathscr{O}_{\mathbb{P}_S^1}).$ The above two functors are pseudo-inverses and hence $\widetilde{\mathscr{D}}_{1,3,d_{1},d_{2}}$ is isomorphic to the scheme $\mathbb{U}$. Note that group sheaf $\textbf{\underline{Aut}}(\mathbb{P}_{\mathbb{Z}}^n, \mathscr{O}(-l_{1}))\times_{PGL_{2}} \textbf{\underline{Aut}}(\mathbb{P}_{\mathbb{Z}}^n, \mathscr{O}(-l_{2}))$ which assigns to every scheme $S$ the group of isomorphisms  $$\phi_{1}: (\mathbb{P}_S^1, \mathscr{O}(-l_{1}))) \to (\mathbb{P}_S^1, \mathscr{O}(-l_{1}))),\phi_{2}:(\mathbb{P}_S^1, \mathscr{O}(-l_{2}))) \to (\mathbb{P}_S^1, \mathscr{O}(-l_{2}))) $$ is isomorphic to $GL_{2}/\mu_{l_{1}}\times_{PGL_{2}} GL_{2}/\mu_{l_2}\cong \Gamma(l_{1},l_{2})$ (see \cite{AV} Theorem $6.15$). Moreover since every pair $(P \to S, \mathscr{L}_{i},i=1,2)$ where $P \to S$ is Brauer--Severi and $\mathscr{L}_{i},i=1,2$ is a line bundle which restricts to degree $-l_{i},i=1,2$ in every geometric fibre is \textrm{\'et}ale locally isomorphic to $(\mathbb{P}_S^1, \mathscr{O}(-l_{i}))),i=1,2$, we have that $\widetilde{\mathscr{D}}_{1,3,d_{1},d_{2}}$ is is a principal bundle over $\widetilde{\mathscr{D}}_{1,3,d_{1},d_{2}}$. Hence ${\widetilde{\mathscr{D}}_{1,3,d_{1},d_{2}}} \cong \mathbb{U}$. Clearly the action is given by $$A \cdot f_{1}(x) = f_{1}(A^{-1} \cdot x),\ (\alpha A) \cdot f_{2}(x) = f_{2}((\alpha A)^{-1}(x)) = (\alpha)^{-l_{2}} \cdot f_{2}(A^{-1}(x)) ,$$ where we identify $\Gamma(l_{1},l_{2})\cong GL_{2}/\mu_{l_{1}}\times _{PGL_{2}} GL_{2}/\mu_{l_{2}}$ by sending
$[(\alpha, A)] \mapsto ([A],[\alpha A])$.
\QEDB

\vspace{5pt}

From the above quotient structure we obtain the following
\begin{theorem}
$\mathscr{D}_{1,3,d_{1},d_{2}}$ has trivial generic stabilizer.
\end{theorem}
\noindent\textbf{Proof.} In order to prove the result we show that the stabilizer of a general geometric point in $\mathbb{U}$ is contained in  $\mu_{l_{1}}\times \mu_{l_{2}}$ where we have the canonical embedding given by
$i:\mu_{l_{1}}\times \mu_{l_{2}} \xhookrightarrow{} G_{m}\times GL_{2}.$ 

Let $(f_{1},f_{2})\in \mathbb{U}$ be a general pair of smooth forms with degree $l_{i}\geq 4, i=1,2$, and
let $(\alpha, A) \in G_{m}\times GL_{2}$ be an element of the stabilizer i.e.
$(\alpha, A)\cdot (f_{1},f_{2})=(f_{1},f_{2})$ which is equivalent to
$$f_{1}(A^{-1}x)=f_{1}(x)\quad\textrm{ and }\quad f_{2}((\alpha. A)^{-1}x)=f_{2}(x)\quad \forall x \in \mathbb{P}^1.$$
We note that a general smooth polynomial of degree $d$ has stabilizer group $\mu_{d}$ and consequently
$A \in \mu_{l_{1}}$ and $\alpha A \in \mu_{l_{2}}$. Thus, we obtain
$$A=\zeta_{l_{1}}^{k_{1}}\quad\textrm{ and }\quad  \alpha=\zeta_{l_{2}}^{k_{2}}/\zeta_{l_{1}}^{k_{1}}$$ where $\zeta_{l_i}$ is a {\it primitive} $l_i$-th root of unity for $i=1,2$.
This immediate implies $[(\alpha, A)]=[1] \in \Gamma(l_{1},l_{2})$
 and hence the action of $\Gamma(l_{1},l_{2})$ is free on some open subset of $\mathbb{U}.$ It follows that the stack $\mathscr{D}_{1,3,d_{1},d_{2}}$ is an orbifold.\QEDB

\subsection{Gerbe Structure of \texorpdfstring{$\mathscr{H}_{1,3,d_{1},d_{2}}$}{} over \texorpdfstring{$\mathscr{D}_{1,3,d_{1},d_{2}}$}{}}\label{H over D}
Consider a natural morphism of stacks $F: \mathscr{H}_{1,3,d_{1},d_{2}} \to \mathscr{D}_{1,3,d_{1},d_{2}}$ which is defined as follows.
\begin{itemize}
    \item[(1)] For an object $\zeta = (P \to S, \mathscr{L}_1, \mathscr{L}_2, i_1, i_2)$ in $\mathscr{H}_{1,3,d_{1},d_{2}}(S)$, $$F(\zeta) = (P \to S, \mathscr{L}_1^{\otimes 2} \otimes \mathscr{L}_2^{\otimes {-1}}, \mathscr{L}_2^{\otimes 2} \otimes \mathscr{L}_1^{\otimes {-1}}, i_1 \otimes \textrm{id}_{\mathscr{L}_2^{\otimes{-1}}}, i_2 \otimes \textrm{id}_{\mathscr{L}_1^{\otimes{-1}}}).$$
    \item[(2)] For a morphism $\Phi = (\lambda_1: P \to P', \lambda_2: \mathscr{L}_1 \to \lambda_1^*\mathscr{L}_1', \lambda_3: \mathscr{L}_2 \to \lambda_1^*\mathscr{L}_2')$, $$F(\Phi) = (\lambda_1: P \to P', \lambda_2^{\otimes 2} \otimes (\lambda_3^*)^{\otimes{-1}}, \lambda_3^{\otimes 2} \otimes (\lambda_2^*)^{\otimes{-1}})$$ 
where $(\lambda_i)^*$ denotes the dual of $\lambda_i$ for $i = 2,3$.

\end{itemize}

\begin{theorem}\label{stacks of cyclic triple covers is a gerbe}

The morphism $F$ realizes $\mathscr{H}_{1,3,d_{1},d_{2}}$ as a gerbe over $\mathscr{D}_{1,3,d_{1},d_{2}}$ with relative automorphism group $\mu_3$. In particular this implies that $\mathscr{D}_{1,3,d_{1},d_{2}}$ is a $\mu_{3}$-rigidification of $\mathscr{H}_{1,3,d_{1},d_{2}}.$

\end{theorem}

\noindent\textbf{Proof.} The proof that $\mathscr{H}_{1,3,d_{1},d_{2}}$ is a gerbe over $\mathscr{D}_{1,3,d_{1},d_{2}}$ is to check the conditions $2(a)$ and $2(b)$ in Lemma $8.11.3$ \cite[\href{https://stacks.math.columbia.edu/tag/06NY}{Tag 06NY}]{stacks-project} as in the proof of Theorem ~\ref{stacks of cyclic covers is a gerbe}. 

 Let us prove 2(a). We have to show that after an \'etale pullback \{$S_{i}\rightarrow S$\}, an object $ y\in \mathscr{D}_{1,3,d_{1},d_{2}}(S_i)$ has a pre-image in $\mathscr{H}_{1,3,d_{1},d_{2}}(S_i)$. So, let us assume $y=(\mathbb{P}^{1}_{S}\to S, \mathscr{O}(-2d_1+d_2),\mathscr{O}(-2d_2+d_1),\lambda_1,\lambda_2),$ where $\lambda_1:\mathscr{O}(-2d_1+d_2)\to \mathscr{O}_{\mathbb{P}^{1}_{S}}$ and similar for $\lambda_2$. Then the pre-image is given by $(\mathbb{P}^{1}_{S}\to S,\mathscr{O}(-d_1),\mathscr{O}(-d_2),\beta_1,\beta_2),$ where $\beta_{1}:\mathscr{O}(-2d_1)\to \mathscr{O}(-d_2)$ is clearly given by $\lambda_1\otimes \textrm{id}_{\mathscr{O}(-d_2)}.$ Similar for $\beta_2.$

Let us prove 2(b).
We have to show that \'etale locally every morphism in $\mathscr{D}_{1,3,d_{1},d_{2}}(S)$ has a pre-image in $\mathscr{H}_{1,3,d_{1},d_{2}}(S).$ A morphism $\lambda$ in $\mathscr{D}_{1,3,d_{1},d_{2}}(S)$ is given by $\lambda=(\lambda_1:\mathbb{P}^{1}_{S}\to \mathbb{P}^{1}_{S},\lambda_2:\mathscr{O}(-2d_1+d_2)\to \mathscr{O}(-2d_1+d_2), \lambda_3:\mathscr{O}(-2d_2+d_1)\to \mathscr{O}(-2d_2+d_1)).$ We have to find a morphism $\beta$ after \'etale pull back such that $F(\beta)=\lambda.$

            Note that $\lambda_2,\lambda_3 \in k^{*}.$ Let $$\beta:=(\beta_1=\lambda_1:\mathbb{P}^{1}_{S}\to \mathbb{P}^{1}_{S}, \beta_2: \mathscr{O}(-2d_1+d_2)\to \mathscr{O}(-2d_1+d_2), \beta_3: \mathscr{O}(-2d_2+d_1)\to \mathscr{O}(-2d_2+d_1)),$$ with $\beta_2,\beta_3 \in k^{*}.$ In order to find a solution we have to solve the following the equation clear from the description of $F,$ given by $\beta_{2}^{\otimes 2}\otimes \beta_{3}=\lambda_2, \beta_{3}^{\otimes 2}\otimes \beta_{2}=\lambda_3.$ On an algebraically closed field the solution is possible.
We show that the relative automorphism group is $\mu_3$. If $F(\Phi) = \textrm{id}$, then $\lambda_1 = \textrm{id}$, $\lambda_2^{\otimes 2} = \lambda_3^*$ and $\lambda_3^{\otimes 2} = \lambda_2^*$. Then, since $\lambda_i^{**} = \lambda_i$ for $i = 2,3 $, we have that $$\textrm{Aut}(\zeta/F(\zeta)) = \left\{(\textrm{id}, \lambda_2, (\lambda_2^{\otimes 2})^*) | \lambda_2^{\otimes 3} = \textrm{id} \right\} = \mu_3.$$

Denote the coarse moduli scheme by $M_{1,3,d_{1},d_{2}}$.  Being an orbifold there exists an open dense schematic substack $\mathscr{D}_{1,3,d_{1},d_{2}}^0$ of $\mathscr{D}_{1,3,d_{1},d_{2}}$ isomorphic to an open subscheme $M^{0}_{1,3,d_1,d_2}\xrightarrow{i}M_{1,3,d_{1},d_{2}}.$ This implies that $\mathscr{D}_{1,3,d_{1},d_{2}}$ admits a section from $M^{0}_{1,3,d_1,d_2}$ which we denote by $p:M^{0}\rightarrow \mathscr{D}_{1,3,d_{1},d_{2}}$ given by, $$(s:P\rightarrow M^{0},D_{1},D_{2})$$ where, $\textrm{deg}(D_{1})=l_{1},\textrm{deg}(D_{2})=l_{2}$ on each geometric fibre of the universal Brauer--Severi scheme, $D_{i}$ for $i=1,2$ being the universal divisors on $(P\rightarrow M^{0})$.\QEDB

\subsection{Degree of line bundles on the universal Brauer--Severi \texorpdfstring{$(P \rightarrow M^0)$}{}}
In order to understand the existence of tautological families on open sets of $M_{1,3,d_1,d_2}$, we need to understand the conditions in which there exist line bundles $\mathscr{L}_{i},i=1,2$ of degree $-d_{1},-d_{2}$ on the universal Brauer--Severi given by $(P\rightarrow M^{0}).$
We observe that there exist $s_{1}:\Gamma(l_{1},l_{2})\rightarrow GL_{2}/\mu_{l_1}$ and $s_{2}:\Gamma(l_{1},l_{2})\rightarrow GL_{2}/\mu_{l_2}$ induced by the canonical maps $[(\alpha,A)]\rightarrow [A]$ where, $(\alpha,A)\in G_{m}\times GL_{2}$.
\begin{theorem}\label{Group Homomorphism}
There exist group homomorphisms $\phi: \Gamma(l_{1},l_{2})\rightarrow \Gamma(d_{1},d_{2})$ over $PGL_{2}.$  
\end{theorem}
\noindent\textbf{Proof.} We can write down
$\phi=(\phi_{1},\phi_{2})$ where,
$$\phi_{i}:\Gamma(l_{1},l_{2})\rightarrow GL_{2}/\mu_{d_{i}},\ i=1,2$$ preserving the base $PGL_{2}$.
Let us consider $\Psi:\Gamma(l_{1},l_{2})\rightarrow GL_{2}/\mu_{r}$ preserving the base $PGL_{2}$ for some $r$. One obtains,
$$\Psi([(\alpha,A)])=[f_{(\alpha,A)}A]$$ for some function $f:G_{m}\times GL_{2}\rightarrow G_{m}$.
We need to check that $\Psi$ is a group homomorphism and is well-defined.

\noindent\underline{\textit{Group homomorphism property of $\Psi$:}}
We choose $[(\alpha,A)],[(\beta,B)]\in
\Gamma(l_{1},l_{2})$ to obtain,
$$\Psi([(\alpha,A)]\cdot [(\beta,B)])=\Psi([(\alpha,A)])\cdot\Psi_{1}([(\beta,B)])\implies \Psi([(\alpha\beta,AB)])=\Psi([(\alpha,A)]).\Psi([(\beta,B)])$$
$$\implies [f_{(\alpha\beta,AB)}\cdot AB]=[f_{(\alpha,A)}.f_{(\beta,B)}\cdot AB]\implies f_{(\alpha\beta,AB)}=f_{(\alpha,A)}\cdot f_{(\beta,B)}\cdot \zeta_{r}^{k},$$ for some $k\in \mathbb{Z},\zeta_{r}$ denoting the \textit{primitive} generator of $r$-th root of unity.

Thus, $\Psi$ induces a group homomorphism $f:G_{m}\times GL_{2}\rightarrow G_{m}/\mu_{r}$, consequently
$ (f)^{r}\in \widehat{G_{m}\times GL_{2}}$.
Thus, we obtain
$$f^{r}=(\textrm{det}^{k_{1}},\textrm{det}^{k_{2}}).$$
Recall that the embedding of $i':\mu_{l_{1}}\times \mu_{l_{2}}\xhookrightarrow{}G_{m}\times GL_{2}$ is given by $(x_1, x_2) \mapsto (x_2/x_1, x_1I_{2\times 2})$.
The generators of $\mu_{l_{1}}\times \mu_{l_{2}}$ are given by $$\left\{\left(\dfrac{1}{\zeta_{l_{1}}},\zeta_{l_{1}}\cdot I_{2\times 2}\right), (\zeta_{l_{2}},I_{2\times 2})\right\}$$ where $\zeta_{l_{i}},i=1,2$ denote the \textit{primitive} generator of $\mu_{l_{i}}$-th root of unity.

Next we find the condition that ensures $\Psi$ is a well defined group homomorphism.

\noindent\underline{\textit{$\Psi$ is well defined:}} 
In order for $\Psi$ to be well-defined we must have,
$$\Psi\left(\left[\left(\dfrac{1}{\zeta_{l_{1}}},\zeta_{l_{1}}\cdot I_{2\times 2}\right)\right]\right)=[1]=\Psi\left(\left[\left(\zeta_{l_{2}},I_{2\times 2}\right)\right]\right)\implies \left[f_{\left(\frac{1}{\zeta_{l_{1}}},\zeta_{l_{1}}\right)}\cdot\zeta_{l_{1}}\right]=[1]=\left[f_{\left(\zeta_{l_{2}},1\right)}\right]$$
$$\implies \left(f_{\left(\frac{1}{\zeta_{l_{1}}},\zeta_{l_{1}}\right)}\cdot\zeta_{l_{1}}\right)^{r}=1=\left(f_{\left(\zeta_{l_{2}},1\right)}\right)^{r}\implies \left(\dfrac{1}{\zeta_{l_{1}}}\right)^{k_{1}}\cdot(\zeta_{l_{1}})^{2k_{2}}\cdot (\zeta_{l_{1}})^{r}=1=(\zeta_{l_{2}})^{k_{1}}$$
$$\implies (\zeta_{l_{1}})^{2k_{2}-k_{1}+r}=1=(\zeta_{l_{2}})^{k_{1}}.$$\\
Thus one obtains the following
$$l_{1}\mid 2k_{2}-k_{1}+r \quad\textrm{ and }\quad l_{2}\mid k_{1}.$$
Thus in order to have $\phi:\Gamma(l_{1},l_{2})\rightarrow \Gamma(d_{1},d_{2})$ over $PGL_{2}$, we need to find $k_{1},k_{2},k'_{1},k'_{2}\in \mathbb{Z}$ such that,
$$l_{1}\mid (2k_{2}-k_{1}+d_{1}),\quad l_{1}\mid (2k'_{2}-k'_{1}+d_{2}),\quad l_{2}\mid k_{1},\quad l_{2}\mid k'_{1}.$$
Now we use the fact that $$l_{1}=2d_{1}-d_{2},\quad l_{2}=2d_{2}-d_{1}.$$
Simplifying further we obtain that solution to the above system is equivalent to correct choices of $n,m$ satisfying the following,
$$n(2d_{1}-d_{2})=(2k_{2}-k_{1}+d_{1}),\quad (2d_{2}-d_{1})\mid k_{1},\quad m(2d_{1}-d_{2})=(2k'_{2}-k'_{1}+d_{2}),\quad (2d_{1}-d_{2})\mid k'_{1}.$$
On further notice one obtains,
$$(2n-1)d_{1}-nd_{2}+k_{1}=2k_{2},\quad 2md_{1}-(m+1)d_{2}+k'_{1}=2k'_{2},\quad (2d_{2}-d_{1})\mid k_{1},\quad (2d_{2}-d_{1})\mid k'_{1}.$$
To find the correct $m,n$ we deal this case by case as follows.

Case (1) $d_{1},d_{2}$ both even: choose $k_{1}=k'_{1}=0,$ then for any choice of $n,m$ we find $k_{2},k'_{2}\in \mathbb{Z}$ satisfying the above conditions.

Case (2) $d_{1},d_{2}$ both odd: choose $k_{1}=k'_{1}=0,$ then choice of $n,m$ both odd gives us $k_{2},k'_{2}\in \mathbb{Z},$ solving the above conditions.

Case (3) $d_{1}$ even, $d_{2}$ odd: choose $k_{1}=k'_{1}=0$, then any even $n$ and odd $m$ gives us $k_{2},k'_{2}\in \mathbb{Z}$ solving the above conditions.

Case (4) $d_{1}$ odd, $d_{2}$ even: choose $k_{1}=l_{2},k'_{1}=0$ then choice of any $n,m$ gives us $k_{2},k'_{2}\in \mathbb{Z}$ solving the above conditions. \QEDB

\begin{remark}\label{zlt13}
We observe that the universal Brauer--Severi $(P\rightarrow M^0)$ is Zariski locally trivial if and only if $\mathscr{O}(1)\in \textrm{Pic}^{\Gamma(l_{1},l_{2})}(\mathbb{P}^1)$, i.e., if and only if there exists a group homomorphism $\Gamma(l_{1},l_{2})\rightarrow GL_{2}.$

We observe from the proof of Theorem ~\ref{Group Homomorphism} that for the homomorphism to exist we can take $r=1.$ Thus the universal Brauer--Severi scheme trivializes if and only if there exist two integers $k_{1},k_{2}\in \mathbb{Z}$ such that $l_{1}\mid (2k_{2}-k_{1}+1)$ and $l_{2}\mid k_{1}$, or equivalently there exist $m,n\in \mathbb{Z}$ satisfying
\begin{equation}\label{eql}
    k_{1}=ml_{2},\quad\textrm{ and }\quad nl_{1}=2k_{2}-ml_{2}+1.
\end{equation}
Simplifying, one obtains that the above system admits a solution if and only if
$$nd_{2}+md_{1}+1$$ is even. Thus \eqref{eql} has no solution if and only if $d_{1},d_{2}$ both are even. In other cases we find that the solution exists.

To summarize, we have non-trivial Brauer--Severi scheme when, $d_{1},d_{2}$ both are even. The Brauer--Severi scheme trivializes in all other cases of $d_{1},d_{2}.$
\end{remark}

In order to compute the relative Picard group $\textrm{Pic}(P|M^{0}):=\textrm{Pic}(P)/s^{*}\textrm{Pic}(M^{0}),$ we introduce similar constructions as in Lemma \ref{relative picard of universal Brauer Severi} as follows.
\begin{lemma}\label{Picard group of Universal Brauer Severi}
Suppose that the unique section $p: M^{0} \to \mathscr{D}_{1,3,d_{1},d_{2}}$ be denoted by the object $(P \to M^0, D_{1},D_{2})$ where $P \to S$ is the universal Brauer--Severi and $D_{i}$ is the universal divisor of degree $l_{i},i=1,2$ on every geometric fibre. Let $U \subseteq M^0$ Then $$\textrm{Pic} (P|_{U}/U) \cong \textrm{Pic}^{\Gamma(l_{1},l_{2})}(\mathbb{P}^1)/K$$
\textcolor{black}{where $K \subseteq \textrm{Lin}^{\Gamma(l_{1},l_{2})}(\mathscr{O}_{\mathbb{P}^1})$ and $\textrm{Lin}^{\Gamma(l_{1},l_{2})}(\mathscr{O}_{\mathbb{P}^1}) = \hat{G}$ is the subgroup of pairs $(\mathscr{O}_{\mathbb{P}^1}, f) $ where $f$ is a $\Gamma(l_{1},l_{2})$ linearization of $\mathscr{O}_{\mathbb{P}^1}$}.
\end{lemma}
\noindent\textbf{Proof.} Let $\mathbb{D}_{i},i=1,2$ be the incidence divisors inside $\mathbb{U} \times \mathbb{P}^1$ defined by $\mathbb{D}_{1}:= \left\{(x_{1},x_{2},y) \in \mathbb{U} \times \mathbb{P}^1 | y \in V_{+}(x_{1}) \right\}$ and similarly one defines $\mathbb{D}_{2}:= \left\{(x_{1},x_{2},y) \in \mathbb{U} \times \mathbb{P}^1 | y \in V_{+}(x_{2}) \right\}$. There is a natural free action of $\Gamma(l_{1},l_{2})$ on $\mathbb{U}$ such that $M^0$ is its geometric quotient. Let $\Gamma(l_{1},l_{2})$ act diagonally on $\mathbb{U} \times \mathbb{P}^1$ so that the action is free and hence admits a geometric quotient $P'$. Now clearly $\mathbb{D}_{i},i=1,2$ are invariant under the action of $\Gamma(l_{1},l_{2})$ and the induced action of $\Gamma(l_{1},l_{2})$ is free. Hence there exists a geometric quotient $D'_{i},i=1,2$. This gives us a section $q: M^0 \rightarrow \mathscr{D}_{1,3,d_{1},d_{2}}$ given by the object $(P'\rightarrow M^0, D_{1}, D_{2}).$ Following the argument of Lemma \ref{relative picard of universal Brauer Severi} the rest of the proof follows similarly.\QEDB

\subsection{Tautological families on Zariski open sets of the coarse moduli of \texorpdfstring{$\mathscr{H}_{1,3,d_{1},d_{2}}$}{}}
In order to lift the unique section $p: M^{0} \to \mathscr{D}_{1,3,d_{1},d_{2}}$ given by the object $(P \to M^0, D_{1},D_{2})$ where $P \to S$ is the universal Brauer--Severi and $D_{i}$ is the universal divisor of degree $l_{i},i=1,2$ on every geometric fibre, to $\mathscr{H}_{1,3,d_{1},d_{2}}$ over some open subscheme of $M^0$ we need to find line bundles $\mathscr{L}_{1},\mathscr{L}_{2}$ of degree $-d_{1},d_{2}$ on the universal Brauer--Severi. Note that in case a section exists with a lifting then there exists $\mathscr{L}_{i},i=1,2$ with degrees $-d_{i},i=1,2$ such that $$D_{1}\in |\mathscr{L}_{1}^{\otimes 2} \otimes \mathscr{L}_{2}^{\otimes -1}|,\quad D_{2}\in |\mathscr{L}_{2}^{\otimes 2} \otimes \mathscr{L}_{1}^{\otimes -1}|$$ with degrees $l_{i},i=1,2.$ We denote by $D^0$ the open schematic substack in $\mathscr{D}_{1,3,d_{1},d_{2}}$ such that the natural course moduli map is an isomorphism to $M^0.$
\begin{theorem}\label{Main Theorem}\label{tautgal}
The following statements hold true.
\begin{enumerate}
    \item There exists a tautological family of cyclic triple covers over some open subscheme of the coarse moduli $M_{1,3,d_1,d_2}$ of $\mathscr{H}_{1,3,d_{1},d_{2}}$.
    \item The family is unique when restricted to the complement of all hyperplane sections in $M_{1,3,d_{1},d_{2}}^0$ in the linear systems of each line bundle of order 3.
    \item In particular, if one of the following holds:
    \begin{itemize}
        \item[(i)] $2 \mid  \textrm{gcd}(d_{1},d_{2})$ and $3 \mid \textrm{gcd}(l_{1}-2,l_{2}-2)$; or
        \item[(ii)] \textcolor{black}{$\textrm{gcd}(d_{2},2)=1$ and $$3(2(l_2(l_1+1)+(l_2-2)(d_2+1))-(2l_2-1)(2d_1+1))\nmid$$ $$\text{gcd}\left(\begin{array}{l} (l_2-1)(4d_2-5d_1(d_2+1)-4d_1^{2})+(l_2-1
)(d_1+2)(4d_1-5d_2),\\ -2(l_2-1)(l_1-1)(2d_1+1)+4(l_1-1)(l_2-1)(d_1+2),\\
2(l_1-1)(4d_2-5d_1(d_2+1)-4d_1^{2})+(l_1-1)(2d_1+1)(-5d_2+4d_1)\end{array}\right),$$}
    \end{itemize} 
    then for those open sets over which a tautological family exists, such a family is unique.
\end{enumerate}
  
\end{theorem}
\noindent\textbf{Proof.} (1). Let the unique section to the natural map $p: D^0 \to M^0$ be denoted by $(\pi: P \to M_{1,3,d_1,d_2}^0, D_{i},i=1,2)$ where $(P \to M_{1,3,d_1,d_2}^0)$ is a 
Brauer--Severi scheme of relative dimension $1$ and $D_{i},i=1,2$ are divisors smooth over $M_{1,3,d_1,d_2}^0$ with relative degree $l_{i},i=1,2$ respectively.
 
For $i = 1,2$, we need to find line bundles $\mathscr{L}_i$ on $\pi|_U: P|_U \to U$ of vertical degree $d_i$ on some open set $U \subseteq M_{1,3,d_1,d_2}^0$, such that $D_i|_U$ is in the linear system of $|\mathscr{L}_{i}^{\otimes 2}\otimes \mathscr{L}_{j}^{\otimes -1}|$. By Theorem \ref{Group Homomorphism}, we have a group homomorphism $\phi: \Gamma(l_1,l_2) \to \Gamma(d_1,d_2)$ preserving $PGL_{2}$ for any $d_{1},d_{2}$ which implies that $\mathscr{O}_{\mathbb{P}^1}(d_i) \in \textrm{Pic}^{\Gamma(l_1,l_2)}(\mathbb{P}^1)$ for $i = 1,2$. By Lemma \ref{Picard group of Universal Brauer Severi} this is equivalent to the existence of line bundles $\mathscr{L}_{i},i=1,2$ with degrees $d_{i}$ on $(P \to M_{1,3,d_1,d_2}^0)$.
Notice that $$\mathscr{L}_{1}^{\otimes 2}\otimes \mathscr{L}_{2}^{\otimes -1}\otimes \mathscr{O}(-D_{1})$$ being a degree 0 line bundle on each geometric fibre on $(P \to M_{1,3,d_1,d_2}^0)$, using Lemma \ref{Picard group of Universal Brauer Severi}, we have $$\mathscr{L}_{1}^{\otimes 2}\otimes \mathscr{L}_{2}^{\otimes -1}\otimes \mathscr{O}(-D_{1})\cong\pi^{*}(\mathscr{M}),$$ where $\mathscr{M}\in \textrm{Pic}(M^0).$ On shrinking $M^0$ further one can assume $D_1$ is in the linear system of $|\mathscr{L}_{1}^{\otimes 2}\otimes \mathscr{L}_{2}^{\otimes -1}|$. Arguing similarly one can assume $D_2$ is in the linear system of $|\mathscr{L}_{2}^{\otimes 2}\otimes \mathscr{L}_{1}^{\otimes -1}|$. \par 
(2). By \cite{GV1}, Lemma $3.2$, part $\textrm{(ii)}$, one has that if a family exists over an open set $U \subset M_{1,3,d_{1},d_{2}}^0$, the number of such families is the number of $3$-torsion elements of $\textrm{Pic}(U)$. Now note that by excision there is a surjection $\textrm{Pic}(M_{1,3,d_{1},d_{2}}^0) \to \textrm{Pic}(U)$ 
where the kernel consists of all hypersurfaces in $M_{1,3,d_{1},d_{2}}^0-U$. \par
(3). We prove this when $2 \mid \textrm{gcd}(d_{1},d_{2})$. The case when $d_1$ is odd follows similarly. Using (\cite{AV}, Theorem 6.6), one obtains using Smith Normal form that the order of  \textrm{Pic}($M_{1,3,d_{1},d_{2}}^0$) is  $\textrm{gcd}(4(l_{1}-1)(l_{2}-1), 2(l_{1})(l_{1}-1), 2(l_{2})(l_{2}-1))$. 
\textcolor{black}{Let us give an explanation of the above fact. Note that |$\textrm{Pic}(M_{1,3,d_{1},d_{2}}^0)|=|\textrm{Pic}(\mathscr{D}_{1,3,d_{1},d_{2}}^0)|=|\textrm{Pic}(\mathscr{D}_{1,3,d_{1},d_{2}})|$ by Proposition \ref{codimension of stacky locus in orbifold for cyclic triple covers}. 
On the other hand $|\textrm{Pic}(\mathscr{D}_{1,3,d_{1},d_{2}})| = |\textrm{Pic}(\mathscr{H}_{1,3,d_{1},d_{2}})|/3$ which follows from 
Lemma ~\ref{injectivity of Pic maps for cyclic}. When $\textrm{gcd}(d_1,d_2)=2$, to compute $|\textrm{Pic}(\mathscr{H}_{1,3,d_{1},d_{2}})|$ it is enough to observe that we have to compute the Smith Normal Form of the matrix given by, 
$$\begin{vmatrix}
2(l_1-1)&-4(l_1-1)&0\\
4(l_2-1)&-2(l_2-1)&0\\
4d_2-5d_1&4d_1-5d_2&0\\
\end{vmatrix},$$ and we observe that $|\textrm{Pic}(\mathscr{H}_{1,3,d_{1},d_{2}})|=\textrm{gcd}(12(l_2-1)(l_1-1),2(l_1-1)(4d_1-5d_2)+4(l_1-1)(4d_2-5d_1),4(l_2-1)(4d_1-5d_2)+2(l_2-1)(4d_2-5d_1)).$} 
Hence if $3 \mid l_{1}-2$ and $3 \mid l_{2}-2$, the only $3$-torsion element in $\textrm{Pic}(M_{1,3,d_{1},d_{2}}^0)$ is the identity.   \QEDB

\subsubsection{\textcolor{black}{Picard group of $M_{1,3,d_1,d_2}^0$} and non-existence of a tautological family on $M_{1,3,d_1,d_2}^0$}

\textcolor{black}{In this subsection we compute the Picard group of $M_{1,3,d_1,d_2}^0$ required for Theorem ~\ref{Main Theorem}, $(3)$ and} show that there does not exist a tautological family on $M_{1,3,d_1,d_2}^0$. 

\begin{lemma}\label{explaining $F$}
The morphism $F$ in Theorem ~\ref{stacks of cyclic triple covers is a gerbe} is induced by the map $$\textrm{id}: \mathbb{A}_{sm}(1,l_1) \times \mathbb{A}_{sm}(1,l_2) \to \mathbb{A}_{sm}(1,l_1) \times \mathbb{A}_{sm}(1,l_2)$$ which is equivariant under the group homomorphism
$\Phi: \Gamma(d_1,d_2) \to \Gamma(l_1,l_2)$ sending $$([A], [\alpha A]) \to \left(\left[\alpha^{-\frac{d_2}{l_1}}A\right], \left[\alpha^{\frac{2d_2}{l_2}} A\right]\right).$$ Treating $\Gamma(d_1,d_2) \cong \mathbb{G}_m \times GL_2/(\mu_{d_1} \times \mu_{d_2})$, we obtain $$\Phi([(\alpha, A)]) = \left[\left(\alpha^{\frac{d_2}{l_1}+\frac{2d_2}{l_2}}, \alpha^{\frac{-d_2}{l_1}}A\right)\right].$$
Note that for any choice of $\alpha^{\frac{d_2}{l_1}}$ and $\alpha^{\frac{2d_2}{l_2}}$, $\left[\left(\alpha^{\frac{d_2}{l_1}+\frac{2d_2}{l_2}}, \alpha^{\frac{-d_2}{l_1}}A\right)\right]$ is well defined in $\Gamma(l_1,l_2)$. 
\end{lemma}
\smallskip

\noindent\textbf{Proof.} It is clear that $F$ is induced by the identity morphism under some $\Phi: \Gamma(d_1,d_2) \to \Gamma(l_1,l_2)$. We prove the description of $\Phi$. Suppose $\Phi([A], [\alpha A]) = ([\lambda_1 A], [\lambda_2 A])$. The condition of equivariance of the identity map gives
\begin{equation*}
 ((\alpha)^{d_2} \cdot f_{1}(A^{-1}x),(\alpha)^{-2d_2} \cdot f_{2}(A^{-1}x)) =  ( \lambda_1^{-l_1} \cdot f_{1}(A^{-1}x),\lambda_2^{-l_2} \cdot f_{2}(A^{-1}x)) 
\end{equation*}
Hence $[\lambda_1] = \left[\alpha^{-\frac{d_2}{l_1}}\right]$ in $GL_2/\mu_{l_1}$ and $[\lambda_2] = \left[\alpha^{\frac{2d_2}{l_2}}\right]$ in $GL_2/\mu_{l_2}$. Hence our assertion follows.\QEDB

\begin{lemma}\label{injectivity of Pic maps for cyclic}
The natural map: $\textrm{Pic}(\mathscr{D}_{1,3,d_1,d_2}) \to \textrm{Pic}(\mathscr{H}_{1,3,d_1,d_2})$ is an injective homomorphism,
\begin{enumerate}
     \item if $\textrm{gcd}(d_1,d_2)=2$ then the index is 3,
     \item if $\textrm{gcd}(d_2,2)=1$ then the index is $2(l_2(l_1+1)+(l_2-2)(d_2+1))-(2l_2-1)(2d_1+1).$
\end{enumerate}

\end{lemma}

\noindent\textbf{Proof.} We have a commutative diagram of exact sequences, 
\[
\begin{tikzcd}
0 \arrow[r] & \langle \Delta_1, \Delta_2, Z \rangle \arrow[r] \arrow[equal]{d}  & \widehat{\Gamma(d_1,d_2)} \arrow[r] & \textrm{Pic}(\mathscr{H}_{1,3,d_1,d_2}) \arrow[r] & 0 \\
0 \arrow[r] & \langle \Delta_1, \Delta_2, Z \rangle \arrow[r]  & \widehat{\Gamma(l_1,l_2)} \arrow[r] \arrow[u, "\widehat{\Phi}"] & \textrm{Pic}(\mathscr{D}_{1,3,d_1,d_2}) \arrow[r] \arrow[u] & 0
\end{tikzcd}
\]
where the morphism $\widehat{\Phi}$ is induced by the morphism $\Phi$ as in Lemma ~\ref{explaining $F$}, $\Delta_i$ is the inverse image of the discriminant divisor in $\mathbb{A}_{sm}(1,l_i)$ inside $\mathbb{A}_{sm}(1,l_1) \times \mathbb{A}_{sm}(1,l_2)$ and $Z$ is the divisor in $\mathbb{A}_{sm}(1,l_1) \times \mathbb{A}_{sm}(1,l_2)$ consisting of pairs of forms with a common zero. This shows that the right vertical map is injective since $\widehat{\Phi}$ is injective (since $\Phi$ is surjective) and its index is same as the index of $\widehat{\Phi}$. \par

Suppose that $l_1$ and $l_2$ are both even. Note that by \cite{AV}, Section 6, $\widehat{\Gamma(l_1,l_2)}$ is generated by $\langle v_1, v_2 \rangle$ such that $v_1 = l_2e_1 + \frac{l_2}{2}e_2$ and $v_2 = \frac{l_1}{2}e_2$ where $e_1([(\alpha, A)]) = \alpha$ and $e_2([(\alpha, A)]) = \textrm{det}(A)$. 
\begin{equation*}
    \widehat{\Phi}(v_1)([(\alpha, A)])  = (v_1 \circ \Phi)([(\alpha, A)]) = v_1\left(\left[\left(\alpha^{\frac{d_2}{l_1}+\frac{2d_2}{l_2}}, \alpha^{\frac{-d_2}{l_1}}A\right)\right]\right)  = \alpha^{\frac{d_2l_2}{l_1}+2d_2}\alpha^{\frac{-d_2l_2}{l_1}} (\textrm{det}(A))^{\frac{l_2}{2}}  = \alpha^{2d_2} (\textrm{det}(A))^{\frac{l_2}{2}}.
\end{equation*}
Note that since both $l_1$ and $l_2$ are even we have that both $d_1$ and $d_2$ are even as well. Hence by \cite{AV}, Section 6, $\widehat{\Gamma(d_1,d_2)}$ is generated by $\langle u_1, u_2 \rangle$ such that $u_1 = d_2e_1 + \frac{d_2}{2}e_2$ and $\textcolor{black}{u_2} = \frac{d_1}{2}e_2$
Suppose $\widehat{\Phi}(v_1) = a_1u_1+a_2u_2$. Then 
\begin{equation*}
    (a_1u_1+a_2u_2)([(\alpha, A)]) = \alpha^{d_2a_1} (\textrm{det}(A))^{\frac{d_2a_1}{2}} (\textrm{det}(A))^{\frac{d_1a_2}{2}} 
     = \alpha^{d_2a_1} (\textrm{det}(A))^{\frac{d_2a_1+d_1a_2}{2}}.
\end{equation*}
Hence we have $2d_2 = d_2a_1 \implies a_1 = 2$ and $d_2a_1+d_1a_2 = l_2 = 2d_2-d_1 \implies a_2 = -1$.
Therefore $\widehat{\Phi}(v_1) = 2u_1-u_2$. Similarly,
\begin{equation*}
    \widehat{\Phi}(v_2)([(\alpha, A)])  = (v_2 \circ \Phi)([(\alpha, A)]) = v_2 \left(\left[\left(\alpha^{\frac{d_2}{l_1}+\frac{2d_2}{l_2}}, \alpha^{\frac{-d_2}{l_1}}A\right)\right]\right)  =  \alpha^{-d_2} (\textrm{det}(A))^{\frac{l_1}{2}} 
\end{equation*}
Suppose $\widehat{\Phi}(v_2) = b_1u_1+b_2u_2$. Then 
\begin{equation*}
     (b_1u_1+b_2u_2)([(\alpha, A)])  = \alpha^{d_2b_1} (\textrm{det}(A))^{\frac{d_2b_1}{2}} (\textrm{det}(A))^{\frac{d_1b_2}{2}}  = \alpha^{d_2b_1} (\textrm{det}(A))^{\frac{d_2b_1+d_1b_2)}{2}}
\end{equation*}
Hence we have $-d_2 = d_2b_1 \implies b_1 = -1$ and $d_2b_1+d_1b_2 = l_1 = 2d_1-d_2 \implies b_2 = 2$.
Therefore $\widehat{\Phi}(v_2) = -u_1+2u_2$. 

Now, since the determinant of the matrix
\[
\begin{vmatrix}
2 & -1 \\
-1 & 2 \\
\end{vmatrix}
\]
is equal to $3$ which is not a unit in $\mathbb{Z}$, we have that the map $\widehat{\Gamma(l_1,l_2)} \to \widehat{\Gamma(d_1,d_2)}$ is strictly injective and hence the map $\textrm{Pic}(\mathscr{D}_{1,3,d_1,d_2}) \to \textrm{Pic}(\mathscr{H}_{1,3,d,d})$ is also strictly injective. \par

Now assume $l_1$ odd which implies $d_2$ is odd. In this case generators of $\widehat{\Gamma(l_1,l_2)}$ are given by $v_1 = l_2e_1+\frac{l_2(l_1+1)}{2}e_2$ and $v_2=l_1e_2$ and after applying the automorphism $\widehat{\Gamma(d_1,d_2)} \to \widehat{\Gamma(d_1,d_2)}$ induced by $[(\alpha,A)] \mapsto [(\alpha^{-1},\alpha A)]$, the generators of $\Gamma(d_1,d_2)$ are given by $u_1 = d_1e_1+\frac{d_1(d_2+1)}{2}e_2$ and $u_2=d_2e_2$. A similar computation again shows that the map $\textrm{Pic}(\mathscr{D}_{1,3,d_1,d_2}) \to \textrm{Pic}(\mathscr{H}_{1,3,d_1,d_2})$ induced by the matrix 
\textcolor{black}{\[
\begin{vmatrix}
2l_2-1 & (l_2(l_1+1)+(l_2-2)(d_2+1))/2 \\
4 & -1-2d_1 \\
\end{vmatrix}
\]}
is strictly injective because the determinant is not an unit and that completes the proof. \QEDB

\vspace{5pt}

\begin{proposition}\label{codimension of stacky locus in orbifold for cyclic triple covers}
The codimension of $M_{1,3,d_1,d_2}-M_{1,3,d_1,d_2}^0$ is greater than or equal to two. Hence $$\textrm{Pic}(M_{1,3,d_1,d_2}^0) = \textrm{Pic}(\mathscr{D}_{1,3,d_1,d_2}^0) = \textrm{Pic}(\mathscr{D}_{1,3,d_1,d_2}).$$ 
\end{proposition}
\noindent\textbf{Proof.} Let $\mathbb{U}^0 \subset \mathbb{U}$ denote the open set of forms $(f_1,f_2)$ with trivial stabilizer under the action of $\Gamma(l_1,l_2)$. It is enough to show that the complement of $\mathbb{U}^0$ has codimension greater than or equal to two in $\mathbb{U}$. Let $U \subset \mathbb{A}_{sm}(1,l_1) \times \mathbb{A}_{sm}(1,l_2)$ be the open set consisting of forms $(f_1,f_2)$ with trivial stabilizer under the action of $\Gamma(l_1,l_2)$. It is enough to show that the complement of $U$ has codimension at least two. This follows from the fact that $U$ contains $\mathbb{A}_{sm}(1,l_1)^0 \times \mathbb{A}_{sm}(1,l_2)^0$, which by Proposition ~\ref{codimension of stacky locus of orbifold for $n = 1$} has a complement of codimension greater than or equal to two inside $\mathbb{A}_{sm}(1,l_1) \times \mathbb{A}_{sm}(1,l_2)$. \textcolor{black}{The assertion on Picard groups then follow from the fact that $\mathscr{M}_{1,3,d_1,d_2}^0 \cong \mathscr{D}_{1,3,d_1,d_2}^0$ and $\mathscr{D}_{1,3,d_1,d_2}$ is smooth.} \QEDB

\vspace{5pt}
   
The following \textcolor{black}{proposition} shows that there does not exist a tautological family over $M_{1,3,d_1,d_2}^0$. 

\begin{proposition}\label{non-existence of global tautological family for cyclic}
There does not exist a tautological \textcolor{black}{family} on $M_{1,3,d_1,d_2}^0$.
\end{proposition}

\color{black}
\noindent\textbf{Proof.} The proof is very similar to \textcolor{black}{Proposition} ~\ref{tautological family in any dimension} and we use the notations of Lemma ~\ref{Picard group of Universal Brauer Severi} and its proof. Recall that we have shown in Lemma ~\ref{Picard group of Universal Brauer Severi} that $\textrm{Pic}(P) = \textrm{Pic}^{\Gamma(l_1,l_2)}(\mathbb{U} \times \mathbb{P}^1)$ where $(P \to M_{1,3,d_1,d_2}^0, D_1, D_2)$ is the tuple consisting of the universal Brauer-Severi scheme and the universal divisors of degree $l_1$ and $l_2$ on each fibre given by the isomorphism $M_{1,3,d_1,d_2}^0 \cong \mathscr{D}_{1,3,d_1,d_2}^0$. If there is a tautological family on $M_{1,3,d_1,d_2}^0$, then there should exist line bundles $\mathscr{L}_1$ and $\mathscr{L}_2$ in $\textrm{Pic}(P)$ such that $D_i \in |\mathscr{L}_i^{\otimes 2} \otimes \mathscr{L}_j^{\otimes -1}|$ in $\textrm{Pic}(P)$, or in other words there must exist $\Gamma(l_1,l_2)$ equivariant line bundles $\mathbb{L}_1$ and $\mathbb{L}_2$ on $\mathbb{U} \times \mathbb{P}^1$ such that for $i =1,2$ the $\Gamma(l_1,l_2)$ invariant incidence divisors $\mathbb{D}_i \in |\mathbb{L}_i^{\otimes 2} \otimes \mathbb{L}_j^{\otimes -1}|$ on $\mathbb{U} \times \mathbb{P}^1$ in $\textrm{Pic}^{\Gamma(l_1,l_2)}(\mathbb{U} \times \mathbb{P}^1)$. Once again to show that this cannot happen we first claim that $$\textrm{Pic}^{\Gamma(l_1,l_2)}(\mathbb{U} \times \mathbb{P}^1)= \textrm{Pic}^{PGL_2}(\mathbb{P}(\mathbb{U}) \times \mathbb{P}^1)$$ We use Lemma $2.3$ of \cite{KS}. Note that there is a (free) diagonal action of $\Gamma(l_1,l_2)$ on both $\mathbb{U} \times \mathbb{P}^1$ and $\mathbb{P}(\mathbb{U}) \times \mathbb{P}^1$. Consider the natural map from $$\mathbb{U} \times \mathbb{P}^1 \to \mathbb{P}(\mathbb{U}) \times \mathbb{P}^1$$ which is identity on $\mathbb{P}^1$. This map is equivariant under the action of $\Gamma(l_1,l_2)$. Now note that $\Gamma(l_1,l_2) = GL_{2}/\mu_{l_1} \times_{PGL_2} GL_{2}/\mu_{l_2}$. Hence there is a closed normal subgroup $\mathbb{G}_m/\mu_{l_1} \times \mathbb{G}_m/\mu_{l_2} \hookrightarrow \Gamma(l_1,l_2)$ as pairs of classes of scalar matrices and its induced action preserves the fibres of the equivariant map $\mathbb{U} \times \mathbb{P}^1 \to \mathbb{P}(\mathbb{U}) \times \mathbb{P}^1$. As in Proposition ~\ref{tautological family in any dimension}, one can show that this action on the fibers is free and transitive. Hence the  
the natural $\Gamma(l_1,l_2) $ equivariant map from $$\mathbb{U} \times \mathbb{P}^1 \to \mathbb{P}(\mathbb{U}) \times \mathbb{P}^1$$ is a $\mathbb{G}_m/\mu_{l_1} \times \mathbb{G}_m/\mu_{l_2}$ torsor.
Now considering that $\displaystyle\frac{\Gamma(l_1,l_2)}{\mathbb{G}_m/\mu_{l_1} \times \mathbb{G}_m/\mu_{l_2}} \cong PGL_{2}$, we have by Lemma $2.3$, \cite{KS} an isomorphism of quotient stacks 
$$[\displaystyle\frac{\mathbb{U} \times \mathbb{P}^1}{\Gamma(l_1,l_2)}] \cong [\displaystyle\frac{\mathbb{P}(\mathbb{U}) \times \mathbb{P}^1}{PGL_{2}}]$$
Therefore comparing their Picard groups we have that $$\textrm{Pic}^{\Gamma(l_1,l_2)}(\mathbb{U} \times \mathbb{P}^1)= \textrm{Pic}^{PGL_2}(\mathbb{P}(\mathbb{U}) \times \mathbb{P}^1)$$
Note that $\mathbb{U} \times \mathbb{P}^1$ is an open subset of $\mathbb{P}(\mathbb{A}_{sm}^0(1,l_1)) \times \mathbb{P}(\mathbb{A}_{sm}^0(1,l_2)) \times \mathbb{P}^1$ and its Picard group is generated by $\mathscr{O}_{\mathbb{U} \times \mathbb{P}^1}(1,0,0)$, $\mathscr{O}_{\mathbb{U} \times \mathbb{P}^1}(0,1,0)$ and $\mathscr{O}_{\mathbb{U} \times \mathbb{P}^1}(0,0,1)$.  Now consider the images of the incidence divisors $\mathbb{D}_i$ (which we still call $\mathbb{D}_i$ by a slight abuse of notation) inside $\textrm{Pic}(\mathbb{P}(\mathbb{U}) \times \mathbb{P}^1) = \textrm{Pic}(\mathbb{P}(\mathbb{U})) \times \textrm{Pic}(\mathbb{P}^1)$. We see that their classes are given by $\mathscr{O}_{\mathbb{U} \times \mathbb{P}^1}(1,0,l_1)$ and $\mathscr{O}_{\mathbb{U} \times \mathbb{P}^1}(0,1,l_2)$. Now suppose that there exist line bundles $\mathbb{L}_1 = \mathscr{O}_{\mathbb{U} \times \mathbb{P}^1}(a,b,c)$ and $\mathbb{L}_2 = \mathscr{O}_{\mathbb{U} \times \mathbb{P}^1}(a',b',c')$ such that $\mathbb{D}_i \in |\mathbb{L}_i^{\otimes 2} \otimes \mathbb{L}_j^{\otimes -1}|$. Comparing the first two coordinates this gives us the equations $2a-a' = 1$ and $2a'-a = 0$. But this does not have integer solutions. Hence we are done.  \QEDB

\begin{proposition}\label{M_{1,3,d_1,d_2}}
$\mathscr{H}_{1,3,d_{1},d_{2}}$ is unirational. The coarse moduli $M_{1,3,d_1,d_2}$ is unirational and is fibred over a rational base by fibres which are birational to the homogeneous space $\frac{GL_2/\mu_{l_1} \times GL_2/\mu_{l_2}}{GL_2/\mu_{l_1} \times_{PGL_2} GL_2/\mu_{l_2}}$. If $\textrm{char}(\mathbb{k}) = 0$, $M_{1,3,d_1,d_2}$ is rationally fibred over a rational base.
\end{proposition}

\noindent\textbf{Proof.} Since, $\mathscr{H}_{1,3,d_{1},d_{2}} \cong [\mathbb{U}/\Gamma(d_{1},d_{2})]$, where $\mathbb{U}$ is an open subscheme of $\mathbb{A}_{sm}(1,l_1)^0 \times \mathbb{A}_{sm}(1,l_2)^0$, we have that $\mathscr{H}_{1,3,d_{1},d_{2}}$ admits a dominant one morphism from $\mathbb{U}$ which is rational. Hence $\mathscr{H}_{1,3,d_{1},d_{2}}$ is unirational. \par 
Since we are dealing with the birational properties of $M_{1,3,d_1,d_2}$, we can replace $\mathbb{U}$ by $\mathbb{A}_{sm}(1,l_1)^0 \times \mathbb{A}_{sm}(1,l_2)^0$. Unirationality of $M_{1,3,d_1,d_2}$ then follows since there is a dominant rational map from $\mathbb{A}_{sm}(1,l_1)^0 \times \mathbb{A}_{sm}(1,l_2)^0 \to M_{1,3,d_1,d_2}$. Consider the natural map
\begin{equation*}
   ( \mathbb{A}_{sm}(1,l_1)^0 \times \mathbb{A}_{sm}(1,l_2)^0)/\Gamma(l_1,l_2) \to \mathbb{A}_{sm}(1,l_1)^0 /(GL_2/\mu_{l_1}) \times \mathbb{A}_{sm}(1,l_2)^0 /(GL_2/\mu_{l_2})
\end{equation*}
sending $[(f_1,f_2)] \to ([f_1],[f_2])$. This map is surjective and the fibres are isomorphic to $\frac{GL_2/\mu_{l_1} \times GL_2/\mu_{l_2}}{GL_2/\mu_{l_1} \times_{PGL_2} GL_2/\mu_{l_2}}$. Now the conclusion follows from part $(2)$ of the next \textcolor{black}{lemma}. \QEDB

\begin{lemma}\label{non-solvability and rationality} Two remarks are in order.
\begin{itemize}
    \item[(1)] The group $GL_2/\mu_{l_1} \times_{PGL_{2}} GL_2/\mu_{l_2}\cong \frac{\mathbb{G}_{m}\times GL_{2}}{\mu_{l_{1}}\times \mu_{l_{2}}}$ is not solvable. Note that solvability of the denominator would guarantee the rationality of the fibres irrespective of the characteristic of the base field (see \cite{CZ17}).
    \item[(2)] The homogeneous space $\frac{GL_2/\mu_{l_1} \times GL_2/\mu_{l_2}}{GL_2/\mu_{l_1} \times_{PGL_2} GL_2/\mu_{l_2}}$ is rational if $\textrm{char}(\mathbb{k}) = 0$. 
\end{itemize}
\end{lemma}
\noindent\textbf{Proof.} (1) Assume, $\frac{\mathbb{G}_{m}\times GL_{2}}{\mu_{l_{1}}\times \mu_{l_{2}}}$ is solvable which implies $\mathbb{G}_{m}\times GL_{2}$ is solvable and implies $GL_{2}$ is solvable and we arrive at a contradiction.

(2) We observe that the dimension of $GL_2/\mu_{l_1} \times_{PGL_{2}} GL_2/\mu_{l_2}$ is the dimension of $PGL_2$ which is 3. The assertion follows from Theorem 4.4 \cite{CZ17}.\QEDB

\bibliographystyle{plain}

\end{document}